\documentclass[12pt,a4paper]{article}
\usepackage[utf8x]{inputenx}
\usepackage{ucs}
\usepackage[english]{babel}
\usepackage{amsmath}
\usepackage{amssymb}
\usepackage{amsfonts}
\usepackage{amsthm}
\usepackage{amscd}
\usepackage{amsbsy}
\usepackage{mathabx}
\usepackage{pgfplots}
\usepackage{fancyhdr}
\usepackage{xspace}
\usepackage{pdfpages}
\usepackage{mathrsfs}
\usepackage{pdfpages}
\usepackage{graphicx}
\usepackage{subfig}
\usepackage{xcolor}
\usepackage{ulem}
\usepackage{authblk}
\numberwithin{equation}{section}

\newtheorem{definition}{Definition}[section]
\newtheorem{proposition}{Proposition}[section]
\newtheorem{theorem}{Theorem}[section]

\theoremstyle{definition}
\newtheorem{remark}{Remark}[section]

\usepackage[bookmarks=true]{hyperref}
\usepackage{bookmark}
\DeclareMathOperator{\spn}{\mathrm{Span}}
\DeclareMathOperator{\arctantwo}{arctan2}
\newcommand*{\diff}{\mathop{}\!\mathrm{d}}
\title{A model of reaching via subriemannian geodesics in Engel-type group}
\date{}
\author[1,2]{C. Mazzetti 
\thanks{caterina.mazzetti2@unibo.it}}
\author[2]{A. Sarti\thanks{alessandro.sarti@ehess.fr}}
\author[1]{G. Citti\thanks{giovanna.citti@unibo.it}}
\affil[1]{\small Department of Mathematics, University of Bologna}
\affil[2]{Centre d'Analyse et de Math\'{e}matique Sociales, Sorbonne Universit\'{e}}
\begin{document}
\maketitle
\tableofcontents
\begin{abstract}
In this paper, we propose a model of arm reaching movements expressed in terms of geodesics in a sub-Riemannian space. We will choose a set of kinematic variables to which motor cortical cells are selective with the purpose of modelling the specific task of reaching. Minimizing trajectories will be recovered as suitable geodesics of the geometric spaces arising from the selective behaviour of M1 neurons. We will then extend this model by including the direction of arm movement. On this set, we will define a suitable sub-Riemannian metric able to provide a geometric interpretation of two-dimensional task-dependent arm reaching movements. 
\end{abstract}

\section{Introduction}

The motor cortex is one of the principal brain areas involved in voluntary movements
, nevertheless the question on how the central nervous system selects one specific trajectory of movement is not fully understood (see \cite{scott2004optimal} as a review). Movement planning and control strategies are indeed not directly measurable, 
yet the observation of certain invariant characteristics has provided many modelling insights on this topic (see also  \cite{graziano2002cortical}, \cite{kalaska2009intention}, \cite{harrison2013towards}, \cite{omrani2017perspectives} for a general analysis of the problem). For example, for two-dimensional arm reaching tasks, Abend et al. \cite{Hatf} and Morasso \cite{morasso1981spatial} found stereotypical patterns of movement based on straight paths and bell-shaped velocity profiles, suggesting that the central command for reaching gestures is formulated in terms of hand trajectories in space.
More generally, E. Todorov \cite{todorov2006optimal} argued that, among all possible movements, the brain selects those that meet appropriate optimality criteria (see also Graziano et al. \cite{graziano2002cortical}, \cite{aflalo2006partial}).
Currently, there is a wide variety of models of arm reaching trajectories based on optimality principles, 
so that movements are selected to minimize a particular cost function (see \cite{hogan1984organizing}, \cite{FH}, \cite{uno1989formation}, \cite{kawato1990trajectory}, \cite{flash2007affine},  \cite{biess2007computational}, \cite{berret2008inactivation} and \cite{flash2001computational} as a review). 
One of the best-known model is the minimum hand jerk criteria, developed by Flash and Hogan \cite{FH}. The cost function to be minimized is the square of the rate of change of hand acceleration integrated over the movement execution time:

\begin{equation}\label{xtre_intro}
\frac{1}{2}\int_{0}^{T} \left(\dddot{x}^{2}+ \dddot{y}^{2}\right)dt,
\end{equation}
where $x$ and $y$ are the time-varying hand positions in a Cartesian coordinate system. 
Finding the minimum of the functional \eqref{xtre_intro} is equivalent to assuming that one of the main goals of reaching tasks is to produce the smoothest possible hand motion. The model produces horizontal arm movements that globally fit well with experimental data and with the invariant patterns found in \cite{Hatf,morasso1981spatial}. Shortly after this article, Uno, Kawato and Suzuki \cite{uno1989formation} proposed the minimum torque-change model, consisting of an objective function given by the square of the rate of change of torque generated by muscles. Here, the cost function depends on the nonlinear dynamics of the musculoskeletal system. In a model of 2007, Biess, Liebermann and Flash \cite{biess2007computational} defined geometric properties (path and posture) for three-dimensional pointing movements in terms of geodesic paths with respect to a kinetic energy in a Riemannian configuration space. In this setting, they were able to separately determine the geometrical and temporal movement features, allowing a unification of previous computational models. 
Although for the following cases the main modelling subject is human locomotion, many phenomenological models have been developed by inferring the cost function from behavioural data (see \cite{arechavaleta2006optimizing}, \cite{arechavaleta2008optimality}, \cite{berret2008inactivation},\cite{bayen2009asymptotic}, \cite{chitour2010analysis}, \cite{jean2010optimal},   \cite{chitour2012optimal}). 
%
The approach followed in these articles is the setting of a nonholonomic control system, whose underlying structure is defined in terms of sub-Riemannian geometry (see F. Jean's book \cite{jean2014control} for a complete overview of sub-Riemannian geometry and its applications to motion planning problems). The authors showed the existence of optimal solutions, applied the Pontryagin maximum principle (\cite{pontryagin1987mathematical}) to the control problem, and finally compared the minimizing trajectories with the experimental data.  
In the present paper, following a procedure similar to \cite{jean2010optimal}, we will deduce an energy functional from neurophysiological data (see \cite{schwartz2007useful} as a review) and 
provide a phenomenological model of reaching arising from the sub-Riemannian geometry we set up. Similar problems have been applied also for visual areas (\cite{boscain2014curve}, \cite{boscain2012optimal}, \cite{boscain2010existence}, \cite{duits2014association}, \cite{franceschiello2019geometrical}): illusory contours and perceived curves elaborated in the cortical areas V1/V2 have been described as geodesics  \cite{petitot2008neurogeometrie}, \cite{citti2006cortical}, \cite{boscain2014curve}. 

Sub-Riemannian geodesics have been deeply studied, 
we mention the work of Beals, Gaveau and Greiner who solved the geodesic problem with explicit formulas for the Heisenberg group \cite{beals2000hamilton}, the works of Sachkov and Moiseev \cite{moiseev2010maxwell} and of Duits et al. \cite{duits2014cuspless} for  geodesics in the group of motions of a plane $SE\left(2\right)$ and within $SE\left(d\right)$, Ardentov and Sachkov \cite{ardentov2011extremal} for geodesics in the Engel group and Bravo-Doddoli and Montgomery \cite{bravo2022geodesics} for geodesics in jet space. 
In particular, it is known that 
there may exist ``abnormal" geodesics that do not satisfy the Hamiltonian system associated with the geodesic variational problem (\cite{Montgomery, montgomery1994abnormal}, \cite{bryant1993rigidity}, \cite{sussmann1996cornucopia}, \cite{bonnard2001role}, \cite{monti2014family}, \cite{chitour2006genericity}). 
Montgomery first provided an example of such abnormal minimizers (\cite{montgomery1994abnormal, montgomery1994singular}). 
For distribution of rank 2, Liu and Sussmann \cite{liu1996shortest} introduced a class of abnormal extremals which are always locally length minimizing (see also \cite{bryant1993rigidity} for the rigidity phenomena of singular curves). Engel manifolds \cite{montgomery1999engel}
are foliated by abnormal geodesics \cite{sussmann1996cornucopia} (see also the works \cite{ardentov2011extremal}, \cite{beschastnyi2018left}, \cite{bravo2022geodesics} for a deepened study of geodesics in Engel group). 


Aim of this paper is to propose a model of arm reaching movements inspired by the minimum-jerk model and by a model of functional architecture of the arm area of primary motor cortex M1 (see our previous work \cite{mazzetti2022functional}). Minimizing trajectories found by Flash and Hogan \cite{FH} will be recovered as suitable geodesics of the geometric space we set up in \cite{mazzetti2022functional} and will recall in section \ref{M1_cells_model}. 
We will follow an approach compatible with \cite{aflalo2006partial}: we will choose a set of kinematic variables to which motor cortical cells are selective with the purpose of modelling the specific task of reaching (without constraints on the dynamics of the musculoskeletal system).  

As a starting point, we will consider the geometry arising from the first mono-dimensional kinematic model of cells selective behaviour. 
Motor cortical cells are selective of time, position, velocity and acceleration of the hand (\cite{ashe1994movement}, \cite{schwartz2007useful}), which will be denoted by $\left(t,x,v,a\right)\in\mathbb{R}^4$. The differential constraints relating the kinematic variables endow very naturally the cortical features space with an Engel structure (see \cite{montgomery1999engel} and \cite{Montgomery}, section 6.2.2). 
Once recalled the geometry of the space, 
we propose as reaching trajectories a geodesics subset of the sub-Riemannian flow, which result to be comparable with the solutions obtained through the minimum-jerk model. 
We will then recall the sub-Riemannian model of M1 functional architecture extended to the coding of two-dimensional movement-related features. 
Cells selective tuning was considered 
for the position variables $\left(x, y\right)\in\mathbb{R}^2$ at time $t\in\mathbb{R}$, the movement direction $\theta\in S^1$, and the velocity and acceleration along $\theta$, denoted by $\left(v, a\right)$.
All of these variables give rise to 
the features space 
\begin{equation}
\mathcal{M}= \mathbb{R}^{3}_{\left(x,y,t\right)} \times S^1_{\theta} \times \mathbb{R}^{2}_{\left(v,a\right)} 
\end{equation}
and their  
differential constraints 
induce to consider the vanishing of the following one-forms
\begin{equation}\label{3}
\omega_{1} = \cos\theta dx + \sin\theta dy - vdt,\quad \omega_{2} = -\sin\theta dx + \cos\theta dy,\quad
\omega_{3} = dv -adt.
\end{equation}
We will denote $D^{\mathcal{M}}$ the horizontal distribution belonging to the intersection of the kernels of the 1-forms \eqref{3}. 
In this setting, 
the group ceases to be nilpotent and its properties cannot be traced back to a Heisenberg-type group. 
Here, we will focus on a subset of horizontal curves of this structure by 
defining
the notion of admissible curves as integral curves of the form 
\begin{equation}
\dot{\gamma} (t) = X_1 + k\left(t\right) X_2 + j\left(t\right) X_3,
\end{equation}
where the vector fields $X_1, X_2, X_3$ 
are the generators of the horizontal distribution $D^\mathcal{M}$.
Since H\"{o}rmander condition is no more guaranteed for admissible curves, we will prove a connectivity property and the existence of a minimum admissible curve joining two arbitrary points of the space.   
%
%
As we have previously mentioned (more details will be recalled in section \ref{sub}), 
in a sub-Riemannian context it is not obvious that each (minimizing) geodesic satisfies the associated geodesics equations on the cotangent space. 
To overcome this issue, we will prove that admissible geodesics are regular (see Definition \ref{def_curva_regolare} from \cite{hsu1992}). To prove this result, we will exploit 
Theorem \ref{sogno_curve_singolari} of \cite{citti2022variational} (see also \cite{citti2021variational} for the characterization of singular curves in graded manifolds).
 The regularity of admissible geodesics implies that they are normal (see Theorem \ref{regolar_then_normal} recalled in section \ref{sub} from \cite{Montgomery}), so that they can be actually found as solutions of the geodesics equations. 
In section \ref{results} we will present a qualitative analysis of admissible geodesics. Through a numerical approximation of the solutions of the 
hamiltonian 
equations, we will show how admissible geodesics allow to represent a wide variety of task-related reaching motions. 

The structure of the paper is the following. In section \ref{sub} we report some of the main properties of sub-Riemannian geodesics that we will then adapt to our model. 
In section \ref{M1_cells_model} we recall the minimum-jerk model and the geometry of M1 functional architecture
on which we base our study. 
In section \ref{math_model} we develop our model of reaching through the analysis of admissible geodesics. 
Section \ref{results} presents part of the results 
where we provide a geometric interpretation of some task-dependent arm reaching movements. Finally, section \ref{concl} summarizes the model proposed in this paper.

%

\section{Sub-riemannian geodesics}\label{sub}
In this section we recall some properties of the sub-Riemannian metric which will be used to express the model (we refer to \cite{jost2008riemannian} and \cite{Montgomery, le2010lecture, agrachev2019comprehensive} for a more detailed presentation). 

\begin{definition}\label{def:distribution}
Let $M$ be a differentiable manifold of dimension $n$. A distribution $D$ is a subbundle of the tangent bundle $TM$, i.e.  
at every point $q\in M$ there exists a neighbourhood $U_{q}\subset M$ and $k$ linearly independent smooth vector fields $X_{1},\cdots ,X_{k}$ defined on $U_{q}$, such that, for any point $p \in U_{q}$
$$Span\left(X_{1_{|p}}, \dots , X_{k_{|p}}\right)= D_{p}\subseteq T_{p}M.$$
The vector space $D_{p}$ is called \textit{horizontal tangent
space} at the point $p$. 
The distribution $D$ defined in this way is called horizontal tangent bundle 
of rank $k$. 
\end{definition}

\begin{definition}
We say that a distribution $D$ is of type $\left(k, n\right)$ if the horizontal distribution has dimension $k$ and its generated Lie algebra has dimension $n$. 
We also recall that a distribution $D$ is bracket generating (or equivalently, that the H\"{o}rmander condition is satisfied) 
if its generated Lie algebra coincides with the tangent space at every point.
\end{definition}

\begin{definition}\label{subRmanifold}
A sub-Riemannian manifold is a triple $(M, D, \langle\cdot,\cdot \rangle_{g})$, where $M$ is a differentiable manifold, $D$ is a
bracket generating distribution 
and $\langle\cdot,\cdot \rangle_{g}$ is a scalar product on $D$. 
\end{definition}

The geometry of the space is described through curves whose tangent vectors belong to the fixed distribution $D$. In particular 
\begin{definition}
 A curve $\gamma:[a,b]\to M$ 
is called horizontal if it is absolutely continuous and $\dot \gamma(t) \in D_ { \gamma (t)}$ for every $t$.  
\end{definition}

Under H\"{o}rmander assumption, a connectivity property holds true 
\begin{theorem}\label{chow}(Chow \cite{chow2002systeme}). If a subbundle $D$ of the tangent bundle of a connected manifold $M$ is bracket generating, then any couple of points can be joined
by a horizontal path.
\end{theorem}
\begin{definition}
The length of a horizontal curve $\gamma:[a,b]\to M$ is defined by
\begin{equation}\label{length_teoria}
l\left(\gamma\right)= \int_{a}^{b} \left\|\dot{\gamma}\left(t\right)\right\|_{g} dt,
\end{equation}
where $\left\|\dot{\gamma}\left(t\right)\right\|_g= \sqrt{\left\langle \dot{\gamma}\left(t\right), \dot{\gamma}\left(t\right)\right\rangle_{g}}$ is computed using the inner product on the horizontal space $D_{\gamma\left(t\right)}$. 
\end{definition}

Thanks to the connectivity condition it is possible to define a distance function as follows:
\begin{equation}\label{carnot-car_distance}
d_c\left(p, q\right):= \inf\left\{l\left(\gamma\right): \gamma\; \text{is a horizontal curve connecting}\: p\;\text{and}\; q\right\}
\end{equation}
This is usually called the Carnot-Carathéodory distance.

\begin{definition}
The horizontal path that realizes the Carnot-Carathéodory distance \eqref{carnot-car_distance} is called geodesic.
\end{definition} 
Montgomery proved the existence of length-minimizers (see Appendix E of \cite{Montgomery}), so that the $\inf$ in \eqref{carnot-car_distance} can be replaced by a minimum. 
\subsection{Geodesics equations}
The sub-Riemannian geodesic flow is governed by a Hamiltonian system (see section 1.5 and Appendix A of \cite{Montgomery} for a detailed description). 
Let us consider a local frame $\left(X_a\right)_{a= 1}^{k}$ of vector fields for the distribution $D$ and its dual on the cotangent bundle, given by
$P_{X_{a}}\left(q,p\right)= p\left(X_a\left(q\right)\right),\: q\in M, p\in T_{q}^{*}M$.

If $g_{ab}\left(q\right)= \left\langle X_a\left(q\right), X_b\left(q\right)\right\rangle_{q}$ is the matrix of inner products defined by the horizontal frame, let us consider $g^{ab}\left(q\right)$ be its inverse $k\times k$ matrix. 

\begin{proposition}\label{prop_hamiltoniana_teorica}
If $P_a$ and $g^{ab}$ are the functions on $T^{*}M$ that are induced by a local horizontal frame $\left(X_a\right)$ as just described, then the Hamiltonian is given by
\begin{equation}\label{hamiltoniana_teorica}
H\left(q,p\right)= \frac{1}{2}\sum g^{ab}P_{a}\left(q,p\right)P_b\left(q,p\right).
\end{equation}
\end{proposition}
Since $X_a= \sum X_a^{i}\left(x\right)\frac{\partial}{\partial x^{i}}$ is the expression for $X_a$ relative to coordinates $x^i$, then $P_{X_{a}}\left(x,p\right)= \sum  X_a^{i}\left(x\right)p_i$, where $p_i:= P_{\frac{\partial}{\partial x^{i}}}$.

\begin{definition}\label{eq_normali}
In terms of the canonical coordinates $\left(x^{i}, p_i\right)\in T^{*}M$, the differential equations governing the geodesics flow, named normal geodesic equations, are given by
\begin{equation}\label{normal_geodesic_equations}
\dot{x}^i= \frac{\partial H}{\partial p_i}\quad,\quad \dot{p}_i= -\frac{\partial H}{\partial x^i}.
\end{equation}
\end{definition}

\begin{definition}
Solutions of \eqref{normal_geodesic_equations} projected on $M$ are called normal geodesics.
\end{definition}

%
Normal geodesics are locally minimizing geodesics, indeed we recall the following
\begin{theorem}(Montgomery \cite{Montgomery}).\label{local_opt_normal_geo}
Let $\left(\gamma\left(s\right), p\left(s\right)\right)$ be a solution of system \eqref{normal_geodesic_equations} on $T^{*}M$. Then every sufficiently short arc of $\gamma$ is a minimizing subriemmannian geodesic. Moreover, $\gamma$ is the unique minimizing geodesic joining its endpoint.
\end{theorem}

\subsection{Regular and singular curves}\label{regular_singular_curver_hsu_gio}
Unlike the Riemannian enivronment, in sub-Riemannian geometry there exist minimizing curves for the length functional \eqref{carnot-car_distance} which are not solutions of the corresponding geodesics equations (see for instance \cite{montgomery1994abnormal}, \cite{Montgomery}, \cite{bryant1993rigidity}, \cite{sussmann1996cornucopia}). 
Below we will adopt the approach developed by L. Hsu \cite{hsu1992} (see also \cite{citti2022variational, giovannardi2020variations}).\\

Let $h$ be a Riemannian metric on the whole tangent bundle $TM$. We complete $X_1,\ldots, X_k$ be a basis of $D_p$, $p\in M$, by adding $X_{k+1},\ldots, X_n$ that generates $\mathcal{V}_p= \left(D_p\right)^{\bot}$.\\
We can therefore express a vector field $V$ in terms of $\left(X_i\right)_{i=1}^n:$
\begin{equation}
V= V_H + V_V = \sum_{i=1}^k v_{H_{i}} X_i + \sum_{j= k+1}^n v_{V_{j}} X_j. 
\end{equation}
With these notations, we express the notion of admissible vector field given by G. Giovannardi \cite{citti2022variational}.
\begin{definition}
Given a curve $\gamma: I\to M$, a vector field $V$ along $\gamma$ with compact support in $I$ is called admissible if it satisfies the following $\left(n-k\right)$ linear first order ordinary differential equations 
\begin{equation}\label{sogno_matrice}
V_{V}'= -BV_{V} - AV_{H}, 
\end{equation}
where $B\left(s\right)$ is a square matrix $\left(n-k\right)\times\left(n-k\right)$ and $A\left(s\right)$ is of order $\left(n-k\right)\times k$, with components
\begin{align*}
a_{ri}= \left\langle \left[\gamma', X_i\right], X_r\right\rangle,\quad
b_{rj}= \left\langle \left[\gamma', X_j\right], X_r\right\rangle,\quad r,j&= k+1, \ldots, n\\
i&= 1,\ldots, k.
\end{align*}
\end{definition}

We now introduce the concept of holonomy map first showed by Hsu \cite{hsu1992} in 1992.
\begin{definition}\label{holonomy_map}
Let $\gamma: I\to M$ be a horizontal curve and $\left[a,b\right]\subset I$. Fixed $V_H\in C^{1}\left(\left(a,b\right),D\right)$ and $V_V\left(a\right)= 0$, let us consider the solution $V_V\left(s\right)$ of ODE \eqref{sogno_matrice}. The holonomy map is defined as
\begin{align*}
H_{\gamma}^{a,b}: C^1\left(\left(a,b\right), D\right)&\rightarrow \mathcal{V}_{\gamma\left(b\right)}\\
V_H&\mapsto V_V\left(b\right).
\end{align*}
\end{definition}
\begin{definition}\label{def_curva_regolare}(Hsu \cite{hsu1992}).
In the above conditions, we say that $\gamma$ restricted to $\left[a,b\right]$ is regular if the holonomy map $H_{\gamma}^{a,b}$ is surjective. If the holonomy map is not surjective, we say that $\gamma$ is singular. 
\end{definition}

G. Giovannardi provided a useful criterion of non-regularity of curves, which consists on the following
\begin{theorem}(\cite{giovannardi2020variations}).\label{sogno_curve_singolari}
The horizontal curve $\gamma$ is singular restricted to $\left[a,b\right]$ if and only if there exists a row vector field $\Lambda\left(s\right)\neq 0$ for all $s\in\left[a,b\right]$ that solves the following system
\begin{equation}\label{sistema_magnifico_singolare}
\begin{cases}
\Lambda'\left(s\right)= \Lambda\left(s\right)B\left(s\right)\\
\Lambda\left(s\right)A\left(s\right)= 0.
\end{cases}
\end{equation}
\end{theorem}


Finally, we recall a theorem which clarifies some relations between the different types of geodesics.

\begin{theorem}\label{regolar_then_normal}(\cite{Montgomery}).
Every regular minimizing curve (i.e. regular geodesic) is normal. 
\end{theorem}

The union between singular and regular geodesics provides the whole set of minimizers (see \cite{Montgomery}, section 5.3). Moreover, as we have just enunciated in the previous theorem, regular minimizing geodesics are normal. Nevertheless, the converse inclusion is false in general, since it may exists normal geodesics which are singular.

\section{Background of phenomenological and neural models}\label{M1_cells_model}
We will report below the minimum jerk model and a neurogeometrical model related to the selective behaviour of primary motor cortical cells, with the purpose of providing clear references for the model we will set up in the next section.  
Hence, in subsection \ref{richiamo_fh} we will briefly recall Flash and Hogan model, and in subsections \ref{jbF} and \ref{cortical_features_structure} we will describe the main geometrical structures arising from primary motor cortex (M1) functional architecture. 


\subsection{The minimum jerk model}\label{richiamo_fh}
\begin{figure}[htbp]
\centering
\includegraphics[scale= 0.5]{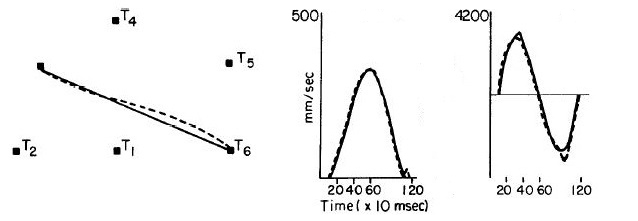}
\caption{Representation example of hand paths,
speed and acceleration for unconstrained point-to-point movement. Dashed lines are the kinematic movements measured. Source: \cite{FH}.}
\label{fm}
\end{figure}

Flash and Hogan assume that movements are planned in terms of hand trajectories rather than joint rotations. Their model is expressed by finding a minimum of an energy function which takes into account the kinematic features of the motion: in moving from an initial to a final position in a given time $T$, the criterion function to be minimized is expressed by

\begin{equation}\label{xtre_2D}
\frac{1}{2}\int_{0}^{T} \left(\dddot{x}^{2}+ \dddot{y}^{2}\right)dt,
\end{equation}
where $x$ and $y$ represent the Cartesian coordinates of hand position. The extremum of the unconstrained cost function, solution of the associated Euler-Lagrange equation,  consists on a fifth order polynomial. 
Assuming that the motion begins and ends with zero velocity and acceleration, the minimum of \eqref{xtre_2D} is given by
\begin{align*}
x\left(t\right)= x_{0} + \left( x_{T}- x_{0}\right)\left( 6 \tau^{5}- 15\tau^{4}+ 10\tau^{3} \right), \quad
y\left(t\right)= y_{0} + \left(y_{T}- y_{0}\right)\left( 6 \tau^{5}- 15\tau^{4}+ 10\tau^{3} \right),
\end{align*}
where $\left(x_0, y_0\right)$ and $\left(x_T, y_T\right)$ are the initial and final hand positions at $t=0$ and $t=T$, and $\tau= \left(\frac{t}{T}\right)$. The model produces straight paths and smooth symmetric velocity profiles that are in accordance with the experimentally observations made by Abend et al. \cite{Hatf} and Morasso \cite{morasso1981spatial} (see Figure \ref{fm} as an example).
The solution trajectories depend only on the initial and
final positions of the hand and movement time, therefore the
optimal trajectory is determined only by the kinematics
of the hand in the task-oriented coordinates and is
independent of the physical system which generates the
motion. 

\subsection{The 1D kinematic tuning model}\label{jbF}
The parameters to which M1 cells are sensible during reaching movements comprise a temporal variable (\cite{georgopoulos1983interruption, kettner1988primate}) together with the speed, acceleration and position of the hand (\cite{ashe1994movement, kettner1988primate, moran1999motor, wang2007motor, truccolo2008primary}). In \cite{mazzetti2022functional} we started by introducing a simple, mono-dimensional model, where the set of kinematic variables can be as
\begin{equation}\label{2_jet_bundle}
\mathcal{J}^{2}=\{\left(t, x\left(t\right), \dot{x}\left(t\right), \ddot{x}\left(t\right)\right)\in\mathbb{R}^{4}| \: t\mapsto x\left(t\right)\in\mathcal{C}^{2}\left(\mathbb{R}\right), t\in\mathbb{R}\}.
\end{equation}
Hence, the resulting structure is described by means of a globally trivial fiber bundle represented by the product $\mathbb{R}^2_{\left(t,x\right)}\times\mathbb{R}^2_{\left(v,a\right)}\simeq \mathcal{J}^2$. The space \eqref{2_jet_bundle} is a two-jet space (see \cite{Montgomery}, section 6.5).   

As remarked in \cite{Encoding},  
the spike probability of a neuron is maximized in the direction of the movement trajectory. In \cite{mazzetti2022functional} we noted that the variables which describe the movement are related by differential constraints and we characterize them by the vanishing of the following 1-forms
\begin{equation}\label{1formaflash}
\begin{cases}
v\diff t- \diff x= 0\\
a\diff t- \diff v= 0.
\end{cases}
\end{equation}
The intersection of the kernels of these 1-forms is a distribution spanned by the vector fields
\begin{equation}\label{vf_x1x2}
X_{1}= \frac{\partial}{\partial{t}}+ v\frac{\partial}{\partial{x}}+ a\frac{\partial}{\partial{v}}, \quad 
X_{2}= \frac{\partial}{\partial{a}}.
\end{equation}
\subsubsection{The Engel group}\label{Engel_structure}
In a neighbourhood of any point of the jet space $\mathcal{J}^2$, the local frame $X_1$, $X_2$ satisfies the property that $X_1$, $X_2$, $\left[X_1, X_2\right]$, $\left[X_1, \left[X_1, X_2\right]\right]$ span the entire tangent bundle. Indeed, from the commutators of the vector fields $X_1$, $X_2$, we get
\begin{align*}
X_{3}:= \left[X_{2}, X_{1}\right] = \frac{\partial}{\partial{v}} , \quad  X_{4} :=\left[X_{3}, X_{1}\right] = \frac{\partial}{\partial{x}} ,\quad 
\left[X_{4}, X_{1}\right] = 
\left[X_{4}, X_{2}\right] = \left[X_{3}, X_{2}\right] = 0.
\end{align*}

These commutation properties define a manifold of Engel type (see \cite{montgomery1999engel} and \cite{Montgomery}, section 6.2.2). 

By choosing on $D$ the metric $g$ which makes $X_1$ and $X_2$ orthonormal, we get a Sub-Riemannian structure on $\mathcal{J}^{2}$. \\
Horizontal curves are integral curves of vector fields $X_1$ and $X_2$, which are of the form 

\begin{equation}\label{geo}
\gamma'\left(s\right)= \alpha_1\left(s\right)X_{1}\left(\gamma\left(s\right)\right)+ \alpha_2\left(s\right)X_{2}\left(\gamma\left(s\right)\right),\\
\end{equation}
where the coefficients $\alpha_i$ are not necessarily constants.\\  
Since the Lie algebra generated by $X_1$ and $X_2$ is the whole tangent space at every point, as recalled in Theorem \ref{chow}, a metric structure is induced on the space. 
%
%

In section \ref{adm_geo_for_center_out} we will express a model inspired by \cite{FH} in the outlined fiber bundle structure, with a choice of 
a horizontal distribution. 

\subsection{The 2D kinematic tuning model of movement directions}\label{cortical_features_structure}
In \cite{mazzetti2022functional} we extended the previous model to a 2D model of the hand movement expressed with the kinematic variables encoded in the brain. Precisely we  considered the 6D space that encodes time, position, direction of movement, velocity and acceleration of the hand in the plane. 
We represented M1 cells tuning variables by the triple $\left(t, x, y\right)\in\mathbb{R}^3$, which accounts for a specific hand's position in time and where the couple $\left(x,y\right)\in\mathbb{R}^2$ represents the cortical tuning for hand's position in a two dimensional space. We also considered the variable $\theta\in S^1$ which encodes hand's movement direction, and the variables $v$ and $a$ which represent hand's speed and acceleration along $\theta$. 
The triple $\left(t, x, y\right)\in\mathbb{R}^3$ is assumed to belong to the base space of the new fiber bundle structure, whereas the variables 
$\left(\theta, v, a\right)\in S^1 \times \mathbb{R}^{2}$ form the selected features on the fiber over the point $\left(t, x, y\right)$. We therefore considered the 6D features set 
\begin{equation}\label{6D}
\mathcal{M}= \mathbb{R}^{3}_{\left(t,x,y\right)} \times S^1_{\theta} \times \mathbb{R}^{2}_{\left(v,a\right)},
\end{equation}

We have then proceeded by exploiting the main differential constraints which characterize the features selected by the single neurons of this area. The angle of movement direction can be defined as $\theta= \arctantwo\left(\dot{y}, \dot{x}\right)$, from which we deduce the equality 
\begin{equation}\label{origine_tutto_0}
\sin\theta\diff x -\cos\theta\diff y= 0.
\end{equation}
Hence, we set the conditions already imposed in section \ref{jbF} 
by means of the vanishing of the 1-forms \eqref{origine_tutto_0} and \eqref{1formaflash}. As a result, we considered the following three 1-forms  
\begin{align}\label{tre_1forme}
\omega_{1} = \cos\theta \diff x + \sin\theta \diff y - v\diff t= 0,\quad
\omega_{2} = -\sin\theta \diff x + \cos\theta \diff y,\quad
\omega_{3} = \diff v -a\diff t
\end{align}
and we searched for the horizontal distribution $D^{\mathcal{M}}$ satisfying the above constraints on the tangent bundle $T\mathcal{M}$. 
The horizontal distribution is given by $D^{\mathcal{M}}= \spn\{X_1, X_2, X_3\}$, where 
\begin{align}\label{campi_2D}
X_{1}= v\cos\theta\frac{\partial}{\partial{x}} + v\sin\theta\frac{\partial}{\partial{y}}+ a\frac{\partial}{\partial{v}}+ \frac{\partial}{\partial{t}},\quad
X_{2}= \frac{\partial}{\partial{\theta}},\quad
X_{3}= \frac{\partial}{\partial{a}}.
\end{align}
Moreover, the following commutation relations 
\begin{equation}\label{commutatori}
\begin{aligned}
\left[X_{1}, X_{2}\right]= v\sin\theta\frac{\partial}{\partial{x}}&- v\cos\theta\frac{\partial}{\partial{y}}=: X_{4},\quad\quad
\left[X_{3}, X_{1}\right]= \frac{\partial}{\partial{v}}=: X_{5},\\
&\left[X_{5}, X_{1}\right]= \cos\theta\frac{\partial}{\partial{x}}+ \sin\theta\frac{\partial}{\partial{y}}=: X_{6},
\end{aligned}
\end{equation}
show that vector fields $\left(X_i\right)_{i=1}^{6}$ are linearly independent. 
Therefore, 
all $\left(X_i\right)_{i=1}^{3}$  together with their commutators span the whole tangent space at every point, meaning that H\"{o}rmander condition holds.\\

We then set on $D^\mathcal{M}$ the metric $g^\mathcal{M}$ which makes $X_1, X_2, X_3$ orthonormal, providing a sub-Riemannian manifold on $\mathcal{M}$.\\

\section{Admissible geodesics for reaching tasks}\label{math_model}
The goal of this section is to develop a model inspired by the phenomenological and neurogeometrical frame recalled in section \ref{M1_cells_model} to describe reaching tasks.
Starting from optimal arm reaching trajectories of Flash and Hogan model \ref{richiamo_fh}, we will lift the problem in the higher dimensional geometric structure introduced in sections \ref{jbF} and \ref{cortical_features_structure}. In this setting, the functional \eqref{xtre_2D} will become an energy functional, whose minima coincide with geodesics. In section \ref{Model}, we will refer to the sub-Riemannian geometry set in section \ref{cortical_features_structure} by analyzing a special case of the associated geodesics problem. These curves will allow a wider variety of reaching movements to be represented. 


\subsection{Kinematic model of 1D motions}\label{model_1D}

In this section, 
we express the functional introduced by Flash and Hogan in terms of the sub-Riemannian jet space introduced in section \ref{jbF}. Minima of functional \eqref{xtre_2D} will coincide with geodesics in this space. 
We will be interested in curves which are lifting of mono-dimensional trajectories. This property will lead to consider only curves in the family of horizontal ones having a non vanishing component along the vector field $X_1$. 
We also refer to \cite{bayen2009asymptotic} for a similar nonholonomic system for modelling human locomotion.

We will then define admissible curves, as follows
\begin{definition}\label{admissible}
A curve $\gamma: \left[0, T\right]\rightarrow \mathcal{J}^2$ is called admissible if it is of the form 
\begin{equation}\label{eq_curve_admissible}
\dot{\gamma} (t) = X_1 + j\left(t\right) X_2.
\end{equation}
\end{definition}
Here, the function $t\mapsto j\left(t\right)$ represents the magnitude of jerk, the rate of change of acceleration.\\
In subsection \ref{Engel_structure} we chose a metric on the distribution $D$ which makes $X_1$ and $X_2$ orthonormal, as a consequence, the length functional on admissible curves reduces to 
\begin{equation}\label{length}
l\left(\gamma\right)= \int_{0}^{T} \sqrt{1+ j^{2}\left(t\right)}\, dt.
\end{equation}
%

We propose here this functional as a good model for voluntary arm reaching movements, since it provides the same solutions as the one presented in \cite{FH}. Indeed, its associated energy functional is
\begin{equation}\label{energy}
E(\gamma)= \frac{1}{2}\int_{0}^{T} \left(1+ j^{2}\left(t\right)\right) dt .
\end{equation}
Since $j(t) = \dddot{x}(t)$ on a lifted curve, there is a strong relation between this functional and the 
one proposed by Flash and Hogan \eqref{xtre_2D}. In order to be able to compare, we will give some geometric properties of the model. 

\bigskip

\begin{remark} \label{commento_singolari_engel}
As outlined in section \ref{regular_singular_curver_hsu_gio}, there could exist minimizers of the length functional which are not solution of the associated hamiltonian system. 
Sussmann in \cite{sussmann1996cornucopia} provided a technique for producing abnormal extremals of arbitrary sub-Riemannian structures on four-dimensional manifolds having a two-dimensional bracket-generating distribution. This result can be applied to our settings, implying that integral curves of the vector field $X_2$ are singular geodesics. Also Byant and Hsu (see \cite{bryant1993rigidity}, Propositions 2.1 and 3.2) explicitly computed their non regularity in the Engel group. 
Moreover, integral curves of $X_2$ are the only singular curves that can be found. Consequently, admissible geodesics on this setting are non singular.
\end{remark}

\subsubsection{Normal geodesics in the two-jets bundle}\label{geo_2jet_structure_general}
Below we will look at some properties of normal geodesics that will be useful in analyzing the admissible ones. 
Since we are in a Lie group, it is not restrictive to consider geodesics starting from the origin, 
and obtain all the others by applying the action of the group.

\begin{proposition}\label{prop1}
A normal geodesic $\gamma$ starting from $0$ is a solution of the following ODE:

\begin{equation}\label{geo_ode_citti}
\gamma' = \Bigg(vp_{x}+ a  (- p_{x}t + p_v(0)) + p_{t} \Bigg) X_1 + \Bigg(p_x \frac{t^2}{2} - p_v(0) t + p_a(0)\Bigg) X_2=
\end{equation}
$$
=p_{t} X_1  + p_a(0) X_2 + 
p_v(0)(a X_1-tX_2) +
p_{x}\left(vX_1 - at X_1 + \frac{t^2}{2} X_2\right), 
$$
for suitable real constants $p_{t} , p_{x}, p_{v}(0), p_{a}(0)$.
\end{proposition}
\begin{proof}
As recalled in Proposition \ref{prop_hamiltoniana_teorica}, normal geodesics solve a ODE system expressed in terms of the 
cotangent coordinates $\left(t, x, v, a, p_{t}, p_{x}, p_{v}, p_{a}\right)\in T^{*}\mathcal{J}^{2}$. 
Since we chose a metric which makes $X_1$, $X_2$ orthonormal,
the Hamiltonian governing the sub-Riemannian geodesic flow on $\mathcal{J}^{2}$ is 
\begin{equation}\label{hamiltoniana}
H= \frac{1}{2}\left((vp_{x}+ ap_{v}+ p_{t})^2 +  p_{a}^{2}\right),
\end{equation}
whose normal geodesic equations are expressed by

\begin{align}\label{ham_flash}
p'_{t}&= 0, \;
p'_{x}= 0, \;
p'_{v}= -p_{x}\left(vp_{x}+ ap_{v}+ p_{t}\right), \;
p'_{a}= -p_{v}\left(vp_{x}+ ap_{v}+ p_{t}\right), \\
t'&= vp_{x}+ ap_{v}+ p_{t}, \;
x'=  v\left(vp_{x}+ ap_{v}+ p_{t}\right), \;
v'= a\left(vp_{x}+ ap_{v}+ p_{t}\right), \;
a'= p_{a}.
\end{align} 
In this way, for the dual variables it holds

$$p_{t}= p_{t}\left(0\right) , \quad p_{x} = p_{x}\left(0\right), \quad p'_{v}= - p_{x}t', \quad p'_{a}= - p_{v}t'= p_{x}t t'- p_v(0) t',$$
so that 
\begin{equation} \label{pav}
{p}_{v} = - p_{x}t + p_v(0), \quad 
{p}_{a} =  p_x \frac{t^2}{2} - p_v(0) t + p_a(0).
\end{equation}

In addition, 
\begin{equation} \label{h}
vp_{x}+ ap_{v}+ p_{t}  = vp_{x}+ a  (- p_{x}t + p_v(0)) + p_{t}.
\end{equation}

As a consequence, the equation satisfied by the observed variables reduces to 

$$
\gamma' = \left(vp_{x}+ ap_{v}+ p_{t}\right) X_1 + p_{a}X_2=$$
$$
\Bigg(vp_{x}+ a  \left(- p_{x}t + p_v(0)\right) + p_{t} \Bigg) X_1 + \left(p_x \frac{t^2}{2} - p_v(0) t + p_a(0)\right) X_2.
$$

\end{proof}

\begin{proposition}
A geodesic $\gamma$ starting from $0$ solution of ODE \eqref{ham_flash} is a horizontal curve. If we call 
$ \alpha_1 $ and $\alpha_2 $ the coefficients of $X_1$ and $X_2$, respectively, we have 
$$\gamma' = \alpha_1 X_1 + \alpha_2 X_2,$$
with 
$\alpha_1^2 + \alpha_2^2 = C,$ where $C$ is a strictly positive constant. 
\end{proposition}
\begin{proof}
Starting from geodesics equation \eqref{geo_ode_citti} we immediately obtain that 
$$\frac{d\alpha_1^2}{ds} = 2 \alpha_1 \alpha_1' = 
2\alpha_1 (v' p_x + a' p_v + a p_v' ) = 2 \alpha_1  p_a p_v = - 2 p_a p_a' = - \frac{d\alpha_2^2}{ds}.$$
\end{proof}

\begin{proposition}\label{p3}
If a geodesic $\gamma$ 
connecting two points $0$ and $\xi_1$ in the interval $[0,T]$ is represented as in the previous 
proposition, then 
$$l(\gamma) = T \sqrt{(\alpha_1^2 + \alpha_2^2)}.$$
\end{proposition}

\subsubsection{Admissible geodesics for center-out movements}\label{adm_geo_for_center_out}

The kinematic properties of one dimensional motions do not depend on the movement direction. Here, we propose how one dimensional movements, which can accordingly reflect a center-out reaching task, are realized by means of admissible geodesics in the 2-jet structure considered.

Therefore, we need to restrict the study of geodesics which are admissible curves. In  particular, it is no more clear if the connectivity property still holds true. 
We will fix an initial value set to the origin, $t(0)=x(0)=v(0)=a(0)=0$, 
arbitrary final values $\left(x_1, v_1, a_1\right)\neq \left(0,0,0\right)$ and we will show the existence of an admissible curve $\gamma$ connecting them. We remark that the time $T$ defined in the functional \eqref{length} is free a priori. However, in the case of connecting two points of the space with an admissible curve, the arrival time is automatically fixed. Therefore, below we will assume to fix the final time and, without loss of generality, set it equal to 1.



%

\begin{proposition}\label{p4}
If we fix the initial value $\xi_0= \left(t_0,x_0,v_0,a_0\right)= \vec{0}$
and an arbitrary final value $\xi_1= \left(1,x_1, v_1, a_1\right)$, then there exist coefficients 
$e_0, e_1, e_2$ such that 

$$\dot{\gamma}(t) = X_1 + (e_0 + e_1 t + e_2 \frac{t^2}{2} )X_2$$
satisfies 
$\gamma(0) =\xi_0$ and  $\gamma(1) =\xi_1.$  
\end{proposition}

\begin{proof}

It is a direct computation that the expression
$$\gamma'(s) = X_1(\gamma(s)) + \left(e_0 + e_1 t(s) + e_2 \frac{t^2(s)}{2} \right)X_2(\gamma(s))$$
implies $t'=1$, so that we can identify the evolution parameter $s$ with the time parameter $t$ and replace the tangent vector $\gamma'$ with the classical $\dot{\gamma}$. For the other components of $\dot{\gamma}$ it holds
\begin{equation}\label{amano}
\begin{cases}
\dot{x}=&  v\\
\dot{v}=& a\\
\dot{a}=&  e_0 + e_1 t + e_2 \frac{t^2}{2}.
\end{cases}
\end{equation}
By integrating \eqref{amano} and imposing the boundary conditions $\gamma\left(0\right)= 0$ and $\gamma\left(1\right)= \xi_1$, the matrix of coefficients $D$ of the linear integrated system \eqref{amano} is invertible, hence proving a direct connectivity result. 

\end{proof}

Consequently, it is possible to define a distance referred to admissible curves: 
\begin{equation}\label{cc_1d_admissible}
d_a\left(\xi_0, \xi_1\right)= \inf\left\{l\left(\gamma\right): \gamma\; \text{is an admissible curve connecting}\: \xi_0\; \text{and}\; \xi_1\right\},
\end{equation}
where $\left(\xi_0,\xi_1\right) = \left(\left(0, 0, 0, 0\right), \left(1, x_{1}, v_{1}, a_{1}\right)\right).$\\
\\
Let us first estimate this distance in terms of the Carnot-Carathéodory distance $d$ as in \eqref{carnot-car_distance}:

%
\begin{proposition}
If $\left(\xi_0, \xi_1\right)= \left(\left(0, 0, 0, 0\right),\left(1, x_{1}, v_{1}, a_{1}\right)\right)$, then
\begin{equation}\label{cc_1d_stima}
d\left(\xi_0, \xi_1\right) \leq d_a\left(\xi_0, \xi_1\right).
\end{equation}
In addition, 
$d_a\left(0, \xi_1\right)\leq k\left\|D^{-1}\right\| |\xi_1|,$
where $k$ is an absolute constant and $D$ is the matrix associated to the integrated system \eqref{amano}.

\end{proposition}
\begin{proof}
The first assertion immediately follows from the definition, whereas for the second one we have
\begin{align*}
d_a\left(0, \xi_1\right) &\leq  \int_0^1 \sqrt{1 + \Big(e_0 + e_1 t + e_2 \frac{t^2}{2} \Big) ^2} dt \leq k\left\|D^{-1}\right\|(1 +  |a_1| + |v_1| +|x_1|).
\end{align*}
\end{proof}
Of course it is not clear if the minimum is achieved or not. However, it exists at any time a (horizontal) curve which is admissible. 

In the following, we will prove that the minimum is attained.

\begin{remark}
A common strategy to show the existence of minimum for the length functional is to study the same problem for the associated energy functional. Indeed, from Montgomery \cite{Montgomery}, it is proved that for a fixed time $T$, length functional minimizers parameterized with constant speed coincide with those of the energy functional. However, in this context we cannot directly apply this proposition, since the set on which the minimum is computed is different (see Definition \ref{cc_1d_admissible}). 
\end{remark}

\begin{proposition}\label{minimal_geo}
In a compact neighbourhood of the origin, there exists a minimal admissible curve (i.e. an admissible geodesic) connecting two points $0$ and 
$\xi_1= \left(1, x_{1}, v_{1}, a_{1}\right).$

\end{proposition}

\begin{proof}
Let us consider the length functional 
$l(\gamma)= \int_{0}^{1} \sqrt{1+ j(t)^{2}}\, dt$
with boundary values $\gamma(0)= 0$, $\gamma(1)= \xi_1$ 
and take a minimizing succession $\gamma_n$ connecting $0$ and $\xi_1$ such that $l\left(\gamma_n\right)\rightarrow\inf l$. As the functional is uniformly bounded, by Ascoli-Arzelà theorem there exists a sub-succession $\gamma_{n_j}$ which uniformly converges to a curve $\gamma$ joining $0$ and $\xi_1$. 
Hence, by the semi-continuity of the length integral it holds $l\left(\gamma\right)\leq \liminf_{j} l\left(\gamma_j\right)$, from which it follows that the minimum is attained.
\end{proof}

\begin{remark}\label{normal_geo_are_adm}
Admissible curves are regular (see Remark \ref{commento_singolari_engel}) and hence normal (see Theorem \ref{regolar_then_normal}), this means that we can search for admissible geodesics through system \eqref{ham_flash}. Moreover, it is possible to explicitly find admissible curves solutions of \eqref{ham_flash} in a neighbourhood of the origin.\\
Let us assume that $t(0)=x(0)=v(0)=a(0)=0$. For every $p_t>0$, there exist a $\delta>0$ and $T>0$ such that, 
for every $p_x, p_v, p_a$ satisfying $|p_x|, |p_v|, |p_a|\leq \delta$, the geodesic found in Proposition \ref{prop1} is an admissible geodesic for every $t\leq T$.

Since we are assuming $v(0)= a(0)=0$ and $p_t$ is a strictly positive constant, the function $h$ defined by $h\left(s\right)= p_{t} + v\left(s\right)p_{x}+ a\left(s\right)p_{v}\left(s\right)$ is different from 0 in a neighbourhood of the origin.
Following the approach used in \cite{HigherElastica} for the study of geodesics in jet spaces, we make a reparameterization of system $\eqref{ham_flash}$ by considering the change of variable
$\frac{\diff}{\diff t}= \frac{1}{h\left(s\right)}\frac{\diff}{\diff s}$,
so that
\begin{equation}\label{duali}
\dot{p}_{a}:= \frac{\diff p_a}{\diff t}= - p_v, \quad \dot{p}_{v} = - p_x,\quad \dot{p}_{x} = 0, \quad \dot{p}_{t} = 0.
\end{equation}
It is therefore immediate that $p_a$ is a polynomial of degree two in the  variable $t$. Moreover, since
$$\dot{h}= \frac{p_a p_v}{h}= - \frac{p_a \dot{p}_a}{h},\quad\text{we get that}\quad \dot{\left(h^2\right)}= -\dot{\left(p_a^2\right)}.$$
In this way, 
$h^2 + p_a^2 = p_t^2 + p_a(0)^2$ in a neighbourhood of the origin. As a consequence, we can express $h$ as a function of $p_a$, for every $t$ such that 
$p_t^2 + p_a(0)^2 -  p_a(t)^2 \geq 0$, 
which means for every $t$ such that 
\begin{equation}\label{dis1}
p_t^2 + p_a(0)^2 -  ( p_x \frac{t^2}{2} - p_v(0) t + p_a(0)) ^2 \geq 0.
\end{equation}
We will choose $T$ as the largest value of $t$ for which \eqref{dis1} is satisfied. 
For such values of $t$, we can 
express $h$ as $h\left(t\right)= \sqrt{p_t^2 + p_a(0)^2  - p_a^2\left(t\right)}$. 
We therefore obtain
\begin{equation}\label{sist_integrato}
\dot{a}:= \frac{\diff a}{\diff t}= \frac{p_a\left(t\right)}{\sqrt{p_t^2 + p_a(0)^2 - p_a^2\left(t\right)}}.
\end{equation}
Through \eqref{duali} and \eqref{sist_integrato} system \eqref{geo} is integrable and the associated solution is an admissible curve. 
\end{remark}

\begin{remark}\label{studio_grafico}
If we fix the same boundary conditions expressed in Flash and Hogan model (see Figure \ref{fm} in \ref{richiamo_fh}), 
from a qualitative study of \eqref{sist_integrato}, 
we can already have a representation of the computed trajectories (see e.g. Figure \ref{Flash_path}). 
Indeed, by considering a movement which starts 
and ends at null velocity and acceleration, since the sign of $\dot{a}$ solely depends on $p_a$ which is a polynomial of degree 2, we get an acceleration profile which has three distinct zeros and one sign change. 
Hence the bell-shaped speed profile is recovered.
\end{remark}  

\begin{proposition}\label{energy_set_length}
The energy functional $E$ defined in \eqref{energy} attains its infimum on the set of minimal admissible geodesics.
\end{proposition}
\begin{proof}
Following the same approach used for example in \cite{Montgomery}, if we call $\sigma$ an admissible 
curve between $0$ and $\xi_1$,  
by Cauchy-Schwarz inequality we get 
\begin{align}\label{csh}
l\left(\sigma\right)^2= \left(\int_0^1 1\cdot\left\|\dot{\sigma}\left(t\right)\right\| dt\right)^2 \leq \int_0^1 \left\|\dot{\sigma}\left(t\right)\right\|^2 dt= 2 E\left(\sigma\right), 
\end{align}
where equality holds 
if and only if $\left\|\dot{\sigma}\left(t\right)\right\|$ is constant. From Proposition \ref{minimal_geo}, we know there exists a minimal admissible geodesic $\gamma$ joining $0$ and $\xi_1$, from which we obtain that
\begin{align*}
E\left(\gamma\right)= \frac{l\left(\gamma\right)^2}{2} \leq \frac{l\left(\sigma\right)^2}{2} \leq E\left(\sigma\right).
\end{align*}
Therefore, for every minimizing admissible geodesic $\gamma$ and any admissible curve $\sigma$, we get $E\left(\sigma\right)\geq E\left(\gamma\right)$. The equality $l\left(\gamma\right)^2= l\left(\sigma\right)^2$ can hold if and only if $\sigma$ is also a geodesic. 
So, unless $\sigma$ is also an admissible geodesic, we have $E\left(\gamma\right)< E\left(\sigma\right)$. 
Finally, since admissible geodesics are normal, it holds  $E\left(\gamma\right)= \frac{l\left(\gamma\right)^2}{2}$
 proving that $E$ really attains the infimum in $\frac{d_a^2}{2}$.
\end{proof}

We have proved that the length and energy functionals have (up to a constant) the same minima also for the family of admissible curves. In this way, we are truly reduced to analyze the functional studied by Flash and Hogan. 

We modelled as admissible geodesics those center-out movements whose kinematic properties are in accordance with the ones experimentally observed in the physical space.
Moreover, singular geodesics would exactly interpret those ``non-admissible'' physical reaching trajectories, since they would be obtained by varying arm's acceleration without changing arm's velocity, nor position. 
%

\subsection{Kinematic model of 2D motions}\label{Model} 
We now refer to the geometrical setting reported in section \ref{cortical_features_structure}, from which we will extend the previous model of reaching by including two-dimensional movement trajectories. 

%

Let us now arrange the Hamiltonian setting which allows to analyze the set of normal geodesics. 

\begin{definition}
Let $\left(X_i\right)_{i=1}^3$ be the vector fields \eqref{campi_2D}. The fiber-linear functions on the cotangent bundle $P_{X_{i}}: T^{*}\mathcal{M}\to \mathbb{R}$ defined by $P_{X_{i}}\left(\eta,p\right)= p\left(X_{i}\left(\eta\right)\right)$ are called the momentum functions for $X_i$.
\end{definition} 
In terms of cotangent coordinates $\left(x,y,t, \theta, v, a, p_{x}, p_{y}, p_t, p_{v}, p_{\theta}, p_{a}\right)\in T^{*}\mathcal{M}$ we can write
\begin{equation}\label{momentum_functions_M}
P_{X_1}= v\cos\theta p_{x} + v\sin\theta p_{y}+ ap_{v}+ p_{t},\quad P_{X_2}= p_{\theta},\quad P_{X_3}= p_{a}.
\end{equation}
\\
Since we have selected on the distribution $D^{\mathcal{M}}$ the metric $g^{\mathcal{M}}$ which makes $\left(X_i\right)_{i=1}^3$ orthonormal, the Hamiltonian function (see Proposition \ref{prop_hamiltoniana_teorica}) simply reduces to the sum of squares of the momentum functions relative to the frame $\left(X_i\right)_{i=1}^3$. 
We therefore enunciate the following
\begin{proposition}(\cite{Montgomery}).
The Hamiltonian governing the sub-Riemannian geodesic flow on $\mathcal{M}$ is 
$H = \frac{1}{2}\left(P_{X_1}^{2}+ P_{X_2}^{2}+ P_{X_3}^{2}\right)$ and
the normal geodesic equations are given by 
\begin{equation}\label{ham_2D}
\begin{cases}
p'_{x}&= 0\\
p'_{y}&= 0\\
p'_{t}&= 0\\
p'_{v}&=  -\left(\cos\theta p_{x}+ \sin\theta p_{y}\right)\left(v\left(\cos\theta p_{x}+ \sin\theta p_{y}\right)+ ap_{v}+ p_{t}\right)\\
p'_{\theta}&=  v\left(\sin\theta p_{x}- \cos\theta p_{y}\right)\left(v\left(\cos\theta p_{x}+ \sin\theta p_{y}\right)+ ap_{v}+ p_{t}\right)\\
p'_{a}&= -p_{v}\left(v\left(\cos\theta p_{x}+ \sin\theta p_{y}\right)+ ap_{v}+ p_{t}\right)\\
x'&=  v\cos\theta\left(v\left(\cos\theta p_{x}+ \sin\theta p_{y}\right)+ ap_{v}+ p_{t}\right)\\
y'&= v\sin\theta\left(v\left(\cos\theta p_{x}+ \sin\theta p_{y}\right)+ ap_{v}+ p_{t}\right)\\
t'&= \left(v\left(\cos\theta p_{x}+ \sin\theta p_{y}\right)+ ap_{v}+ p_{t}\right)\\
v'&=  a\left(v\left(\cos\theta p_{x}+ \sin\theta p_{y}\right)+ ap_{v}+ p_{t}\right)\\
\theta '&= p_{\theta}\\
a'&= p_{a}.
\end{cases}
\end{equation}
\end{proposition}

Our purpose will be to identify suitable subsets of normal geodesics in order to provide a phenomenological description for some relevant task-reaching movements.

\subsubsection{Admissible geodesics: reaching targets with prescribed directions}\label{sub_fixed_directions}
  
In this section we will model a cognitive reaching task in which it is required to grasp a target in a specific hand orientation, knowing the initial hand configuration.
We claim that the lifted curves of the space could represent the  
arm reaching trajectories,
and, as we did for center-out reaching movements, we will consider integral curves of the horizontal distribution with a non vanishing component along the vector field $X_1$.\\  

We will then look for admissible curves, as follows
\begin{definition}\label{admissible_2D}
A curve $\gamma: \left[0, 1\right]\rightarrow\mathcal{M}$ is called admissible if it is of the form 
\begin{equation}\label{gamma_adm_2d}
\dot{\gamma} (t) = X_1 + k\left(t\right) X_2 + j\left(t\right) X_3.
\end{equation}
\end{definition}
Here, the function $t\mapsto k\left(t\right)$ represents the
Euclidean curvature over the path $\left(x,y\right)$, whereas the function $t\mapsto j\left(t\right)$ describes the rate of change of acceleration.\\

We then search for admissible curves joining arbitrary couples of points in the cortical feature space $\mathcal{M}$.
\begin{proposition}\label{prop_connection_curve_amm_2d}
If we fix a 
constant $k\in\mathbb{R}$ and we consider arbitrary values $\left(\eta_0, \eta_1\right)= \left(\left(0, x_0, y_0, \theta_0, v_0, a_0\right), \left(1, x_1, y_1, \theta_1, v_1, a_1\right)\right)\in\mathcal{M}$, then there exist constants 
$j_0, j_1, j_2, j_3$ such that 

\begin{equation}\label{connection_curve_amm_2d}
\dot{\gamma}(t)= X_{1}+ k X_{2} + \left(j_0 + t j_1 + j_2 \frac{t^2}{2}  + j_3 \frac{t^3}{3!}\right)X_3
\end{equation}
satisfies 
$\gamma(0) = \eta_0$ and  $\gamma(1) = \eta_1.$  
\end{proposition}
\begin{proof}
Analogously to Proposition \ref{p4}, equation \eqref{connection_curve_amm_2d} sets up system 
\begin{equation}\label{adm_system_2D}
\begin{cases}
\dot{x}& = v\left(t\right)\cos\left(\theta\left(t\right)\right)\\
\dot{y}& = v\left(t\right)\sin\left(\theta\left(t\right)\right)\\
\dot{\theta}& = k\\
\dot{v}& = a\left(t\right)\\
\dot{a}& = j_0 + j_1 t + j_2 \frac{t^2}{2}  + j_3 \frac{t^3}{3!}
\end{cases}
\end{equation}
which can be explicitly integrated. By imposing initial and final conditions for the joining of points $\eta_0$, $\eta_1$, the integrated equations give rise to a linear system in the variables $\left(j_i\right)_{i=0}^{3}$. It is a direct computation to verify that the matrix associated to the integrated system is invertible and hence to prove the existence of coefficients $\left(j_i\right)_{i=0}^{3}$ for equation of \eqref{connection_curve_amm_2d}.
\end{proof}

Thanks to the connectivity property above exposed, it is possible to define a distance referred to admissible curves in the connected space $\mathcal{M}$:  
\begin{equation}\label{cc_2D_admissible}
d^\mathcal{M}_{a}\left(\eta_0, \eta_1\right)= \inf\left\{l\left(\gamma\right): \gamma\; \text{is an admissible curve connecting}\: \eta_0\; \text{and}\; \eta_1\right\},
\end{equation}
where $\left(\eta_0,\eta_1\right) = \left(\left(0, x_{0}, y_0, \theta_0, v_{0}, a_{0}\right), \left(1, x_{1}, y_1, \theta_1, v_{1}, a_{1}\right)\right)\in\mathcal{M}$ and the length $l$ is given by
\begin{equation}\label{length_a_2d}
l\left(\gamma\right)= \int_{0}^{1} \sqrt{1+ k^2\left(t\right)+ j^{2}\left(t\right)}\, dt,
\end{equation}
with $\dot{\gamma}$ solution of \eqref{gamma_adm_2d}. 

For reader convenience, we outline that the energy functional on admissible curves in $\mathcal{M}$ reduces to
\begin{equation}\label{energy_a_2d}
E(\gamma)= \frac{1}{2}\int_{0}^{1} \left(1+ k^2\left(t\right)+ j^{2}\left(t\right)\right) dt .
\end{equation}


Now we will show that admissible curves can be found as solutions of system \eqref{ham_2D}. To do so, we will prove that admissible curves are regular, in the sense of Definition \ref{def_curva_regolare}.   

As recalled in section \ref{regular_singular_curver_hsu_gio}, L. Hsu gave a characterization for a curve to be regular by means of the holonomy map (see Definition \ref{holonomy_map}) and G. Giovannardi 
proved a criterion for identifying singular curves (Theorem \ref{sogno_curve_singolari}).
By applying these results to our case, we verify that even if we considered the whole sub-Riemannian setting, integral curves of $X_2$ or $X_3$ (those whose directions point along the fiber) are singular (see Remark \ref{remark_singular_1} and \ref{remark_singular_2} for the computations). 
We prove the regularity of admissible curves in Remark \ref{remark_regularity}. 
\begin{remark}\label{remark_singular_1}
Let $\gamma: \left[0,1\right]\rightarrow \mathcal{M}$ be an integral curve of the vector field $X_3$. Then $\gamma$ is singular.
\begin{proof} After splitting $X_3$ along $\gamma$ in its horizontal and vertical part 
\begin{equation}
X_{3_{H}}= \sum\limits_{i=1}^{3} v_{H_{i}} X_i\quad,\quad X_{3_{V}}= \sum\limits_{j=4}^{6} v_{V_{j}} X_j,
\end{equation}
a direct computation shows that the admissibility system expressed through the matrix form \eqref{sogno_matrice} is given by
\begin{equation}
V_{V}'\left(s\right)= -A\left(s\right)
\begin{pmatrix}
v_{H_1}\left(s\right) \\
v_{H_2}\left(s\right) \\
v_{H_3}\left(s\right)\\
\end{pmatrix},
\quad\text{where}\quad
A= 
\begin{pmatrix}
0 & 0 & 0 \\
1 & 0 & 0 \\
0 & 0 & 0 
\end{pmatrix}
\end{equation}
for $v_{H_1}, v_{H_2}, v_{H_3}\in C_{0}^{1}\left(\left(0,1\right)\right)$. Then the image of the holonomy map (\ref{holonomy_map}) is equal to 
\begin{equation}\label{prima_hol_map}
V_{V}\left(1\right)= 
\begin{pmatrix}
0 \\
-\int_0^1 v_{H_1}\left(s\right) ds \\
0\\
\end{pmatrix},
\end{equation}
from which we deduce that the holonomy map is not surjective. Hence $\gamma$ is singular. 
\end{proof}
\end{remark}

\begin{remark}\label{remark_singular_2}
A horizontal curve $\gamma:\left[0,1\right]\rightarrow\mathcal{M}$ solution of 
$\gamma'=k\left(s\right)X_2 + j\left(s\right)X_3$, where $k$ and $j$ are different from zero, is singular. 
\begin{proof}
As before, a straightforward computation reveals that the matrices $A$ and $B$ of the admissibility system \eqref{sogno_matrice} are given by 
\begin{equation}\label{ap0}
A= 
\begin{pmatrix}
-k\left(s\right) & 0 & 0 \\
j\left(s\right) & 0 & 0 \\
0 & 0 & 0 
\end{pmatrix}\quad
B= 
\begin{pmatrix}
0 & 0 & -\frac{k\left(s\right)}{v\left(s\right)} \\
0 & 0 & 0 \\
k\left(s\right)v\left(s\right) & 0 & 0 
\end{pmatrix}.
\end{equation}
To prove the singularity of $\gamma$, we will apply Theorem \ref{sogno_curve_singolari}. We will verify if there exists a row vector field $\Lambda\left(s\right)\neq 0$ for all $s\in\left[0,1\right]$ that solves system \eqref{sistema_magnifico_singolare}.
Since it must hold $\Lambda A= 0$, if $\Lambda$ is of the form $\Lambda\left(s\right)= \left(\lambda_1\left(s\right), \lambda_2\left(s\right), \lambda_3\left(s\right)\right)$, then 
\begin{equation}\label{ap1}
-k\left(s\right)\lambda_1\left(s\right)+ j\left(s\right)\lambda_2\left(s\right)= 0.
\end{equation}
Hence, 
\begin{equation}\label{ap2}
\Lambda B=  \begin{pmatrix}
\lambda_1, & \lambda_2, & \lambda_3
\end{pmatrix}
\begin{pmatrix}
0 & 0 & -\frac{k\left(s\right)}{v\left(s\right)} \\
0 & 0 & 0 \\
k\left(s\right) v\left(s\right) & 0 & 0 
\end{pmatrix}=
\begin{pmatrix}
\lambda_3 k\left(s\right) v\left(s\right), & 0, & -\lambda_1\frac{k\left(s\right)}{v\left(s\right)}
\end{pmatrix},
\end{equation}
which implies that 
\begin{equation}\label{ap3}
\begin{cases}
\lambda_1'\left(s\right) &= \lambda_3\left(s\right) k\left(s\right) v\left(s\right)\\
\lambda_2'\left(s\right) &= 0\\
\lambda_3'\left(s\right) &= -\lambda_1\left(s\right) \frac{k\left(s\right)}{v\left(s\right)}.
\end{cases}
\end{equation}
Consequently, $\lambda_2$ is constant and from \eqref{ap1} and \eqref{ap3} we can integrate $\lambda_3$. 
Therefore, for any choice of $\lambda_2\neq 0$, we find a row vector $\Lambda$ solution of \eqref{sistema_magnifico_singolare} whose components are not null. The curve $\gamma$ is therefore singular. 
\end{proof}
\end{remark}

\begin{remark}\label{remark_regularity}
The admissible curve $\gamma:\left[0,1\right]\rightarrow\mathcal{M}$ solution of $\gamma'= X_1 + k\left(s\right)X_2 + j\left(s\right)X_3$ is regular. 
\begin{proof}
In this case, the matrices of system \eqref{sogno_matrice} are
\begin{align*}
A= 
\begin{pmatrix}
-k\left(s\right) & 1 & 0 \\
j\left(s\right) & 0 & 1 \\
0 & 0 & 0 
\end{pmatrix}\quad
B= 
\begin{pmatrix}
\frac{a\left(s\right)}{v\left(s\right)} & 0 & -\frac{k\left(s\right)}{v\left(s\right)} \\
0 & 0 & 0 \\
k\left(s\right)v\left(s\right) & -1 & 0 
\end{pmatrix}.
\end{align*}
Then, as in the previous remark, we look for a row vector field $\Lambda\neq 0$ which solves system \eqref{sistema_magnifico_singolare}.
Since it must be $\Lambda A= 0$, we have that $\Lambda$ is of the form $\Lambda\left(s\right)= \left(0, 0, \lambda\left(s\right)\right)$.\\
Hence, 
\begin{align*}
\Lambda B=  \begin{pmatrix}
0, & 0, & \lambda\left(s\right)
\end{pmatrix}\begin{pmatrix}
\frac{a\left(s\right)}{v\left(s\right)} & 0 & -\frac{k\left(s\right)}{v\left(s\right)} \\
0 & 0 & 0 \\
k\left(s\right)v\left(s\right) & -1 & 0 
\end{pmatrix}= \begin{pmatrix}
\lambda\left(s\right) k\left(s\right)v\left(s\right), & -\lambda\left(s\right), & 0
\end{pmatrix}.
\end{align*}
This means that $\lambda'\left(s\right)= 0$ and $\lambda\left(s\right)= 0$, therefore the unique solution to system \eqref{sistema_magnifico_singolare} is $\Lambda \equiv 0$. This enables to conclude that admissible curves are regular.  
\end{proof}
\end{remark}

Through Remark \ref{remark_regularity} we have proved that admissible curves are regular, therefore solutions of the hamiltonian system \eqref{ham_2D}, as stated by Theorem \ref{regolar_then_normal}. 
As we did in Remark \ref{normal_geo_are_adm} in section \ref{adm_geo_for_center_out}, we can explicitly represent solutions of \eqref{ham_flash} which are admissible in a neighbourhood of the origin (see Remark \ref{normal_geo_2D_are_adm}).

\begin{remark}\label{normal_geo_2D_are_adm}
Let us assume that $\eta_0= \vec{0}\in\mathcal{M}$. For every $p_t>0$, there exist a constant $k$, a $\delta>0$ and $T>0$ such that, 
for every $p_x, p_y, p_\theta, p_v, p_a$ satisfying $|p_x|, |p_y|, |p_\theta|, |p_v|, |p_a|\leq \delta$, the solution of system \eqref{ham_2D} is an admissible geodesic for every $t\leq T$.
We define the function $\psi\left(s\right)= v\left(s\right)\left(\cos\left(\theta\left(s\right)\right) p_{x}+ \sin\left(\theta\left(s\right)\right) p_{y}\right)+ a\left(s\right)p_{v}\left(s\right)+ p_{t}$ and we parameterize the equations with respect to $t$ by setting $\frac{\diff }{\diff t}= \frac{1}{\psi\left(s\right)}\frac{\diff}{\diff s}.$ 
Since we assumed $v\left(0\right)= a\left(0\right)= \theta\left(0\right)= 0$ and $p_t>0$, the function $\psi$ is strictly positive and $p_\theta\sim k$ in a neighbourhood of the origin.
Hence, the function $p_a$ is a polynomial of degree two in the  variable $t$ in a neighbourhood of $\eta_0$. Moreover, since
$\dot{\psi}=  -\frac{1}{\psi}\left(p_{\theta}\dot{p_{\theta}}+ p_a \dot{p_a}\right)$,
we obtain that $\psi^2 + p_{\theta}^2+ p_a^2 = p_t^2+ k^2+ p_{a}\left(0\right)$. 
We can therefore express $\psi$ as a function of $p_\theta, p_a$ for every $t$ such that 
\begin{equation}\label{dis2}
p_t^2+ k^2+ p_a(0)^2- p_\theta(t)^2 -  p_a(t)^2 \geq 0.
\end{equation} 
By choosing $T$ as the largest value of $t$ for which \ref{dis2} is satisfied, we can 
express $\psi$ as $\psi\left(t\right)= \sqrt{p_t^2+ k^2+ p_a(0)^2- p_\theta(t)^2 -  p_a(t)^2}$, for every $t\leq T$. 
We therefore obtain
\begin{align*}
\dot{\theta}\left(t\right)= \frac{p_{\theta}\left(t\right)}{\sqrt{p_t^2+ k^2+ p_a(0)^2- p_\theta(t)^2 -  p_a(t)^2}},\,
\dot{a}\left(t\right)= \frac{p_a\left(t\right)}{\sqrt{p_t^2+ k^2+ p_a(0)^2- p_\theta(t)^2 -  p_a(t)^2}}.
\end{align*}
\end{remark}

\medskip
Admissible curves in the fiber bundle structure are those that allow to move from one fiber to another. We point out how the choice of variables is based on neurophysiological and physiological findings. Indeed, $\theta$ and $a$ are the variables engrafted in the motor cortex, while the kinematic variables describe movement in the external world space. 

\begin{remark}\label{remark_per_non_ripetere}
The same techniques used in \ref{adm_geo_for_center_out} (see Proposition \ref{minimal_geo} and \ref{energy_set_length}) can be applied to prove that the inf in \eqref{cc_2D_admissible} is a minimum and that minimizing sets for the energy and length functional coincide. In this way, it is still possible to consider admissible geodesics for the space $\mathcal{M}$.
\end{remark}

\subsubsection{Geodesics between sets
}\label{sub_free_direction}
In this section, we will analyze a more general situation. Indeed, we would like to model the circumstance where the object to be reached does not require a particular orientation with which to be grasped, or it is indifferent how to grasp it in terms of acceleration. 
In this case we will impose that the second extreme of the geodesic belong to a set. As a result, the movement trajectory will be defined as the minimizing geodesic between two a priori known sets, obtained by fixing the $\left(x, y, t, v\right)\in\mathbb{R}^{4}$ components and by varying the directions and accelerations $\left(\theta, a\right)\in S^1\times\mathbb{R}$ variables.
This method is the same adopted by B. Franceschiello in \cite{franceschiello2019geometrical} for the modeling of illusory contours in the visual system. 
Below we will show the main definitions and statements.
\begin{definition}\label{definizione_distanza_p_insieme}
Let $F_0\subset\mathcal{M}$ be a compact and non empty set. 
We define the distance function from $F_0$ as
\begin{equation}\label{def_distanza_p_insieme}
d^{\mathcal{M}, F_0}_{a}\left(\eta\right)= \inf_{\eta_0\in F_0} d_a^{\mathcal{M}}\left(\eta_0, \eta\right),
\end{equation}
where $d_a^{\mathcal{M}}\left(\eta_0, \eta\right)$ is the distance referred to admissible curves (see \ref{cc_2D_admissible}) in the cortical space $\mathcal{M}$.
\end{definition}

\begin{definition}\label{definizione_distanza_insieme_insieme}
Let $F_0, F_1\subset\mathcal{M}$ be compact and non empty sets. We define the distance function between $F_0, F_1$ as 
\begin{equation}\label{def_distanza_insieme_insieme}
d^{\mathcal{M}}_{a}\left(F_0, F_1\right)= \inf_{\eta_{1}\in F_1} d^{\mathcal{M}, F_0}_{a}\left(\eta_1\right),
\end{equation}
where $d^{\mathcal{M}, F_0}_{a}$ is the distance function from $F_0$ defined in \eqref{def_distanza_p_insieme}.
\end{definition}

\begin{remark}\label{oss_minimi}
It is clear that if $F_0$ and $F_1$ respectively reduce to $\left\{\eta_0\right\}$ and $\left\{\eta_1\right\}$, then distance $d^{\mathcal{M}}_{a}\left(F_0, F_1\right)$ turns out to be $d_a^{\mathcal{M}}\left(\eta_0, \eta_1\right)$, for which we outline the existence of a minimum in Remark \ref{remark_per_non_ripetere}.
\end{remark}

\begin{definition}\label{def_geo_amm_set}
In the same conditions of Definition \ref{definizione_distanza_insieme_insieme}, we call admissible geodesic with extrema in the sets $F_0$ and $F_1$ the admissible curve $\gamma: \left[0,1\right]\to \mathcal{M}$ such that  $\gamma\left(0\right)\in F_0$, $\gamma\left(1\right)\in F_1$ and which realizes the minimum in \eqref{def_distanza_insieme_insieme}.
\end{definition}

Due to the compactness of $F_0$ and $F_1$, it is immediate to prove the following

\begin{proposition}
In the same conditions of Definition \ref{definizione_distanza_insieme_insieme}, there exists an admissible curve with extrema in $F_0$ and $F_1$ for which the minimum in \eqref{def_distanza_insieme_insieme} is attained
\begin{proof}
We can find two sequences, respectively $\left(\eta_0\right)_n$ in $F_0$ and $\left(\eta_1\right)_n$ in $F_1$, 
such that $d^{\mathcal{M}}_{a}\left(\left(\eta_0\right)_n, \left(\eta_1\right)_n\right)$ tends to $d^{\mathcal{M}}_{a}\left(F_0, F_1\right)$. Since $\left(\eta_0\right)_n$ and $\left(\eta_1\right)_n$ are bounded in a compact set, 
there exist two sub-successions $\left(\eta_0\right)_{n_{j}}$, $\left(\eta_1\right)_{n_{j}}$ in $F_0$ and $F_1$ which uniformly converges to $\eta_0$ and $\eta_1$, respectively. A geodesic between $\left(\eta_0, \eta_1\right)$ exists, as recalled in Remark \ref{oss_minimi}, and attains its minimum in \eqref{def_distanza_insieme_insieme}.
\end{proof}
\end{proposition}
\section{Results}\label{results}
This section is dedicated to some experimental simulations for the solution of systems \eqref{ham_flash} and \eqref{ham_2D}. 
Our goal is to provide a neurogeometrical interpretation of some task-dependent arm reaching movements using properties of geodesics established in \ref{adm_geo_for_center_out}, \ref{sub_fixed_directions} and \ref{sub_free_direction}. 
For each of the cases we analyze, we will assume a fixed initial and final position, together with a null velocity at the beginning and at the end of the movement.
First of all, in section \ref{comp_FH}, we recognize that solutions of  \eqref{ham_flash} and \eqref{ham_2D} projected on the 2D plane coincide with are comparable to the minimizers computed through Flash and Hogan functional introduced on the basis of experimental evidence (see \cite{FH} and \cite{Hatf, morasso1981spatial}).
In order to fall in the assumptions adopted in Flash and Hogan model, we will impose to the equations for the 2D case \eqref{ham_2D} the constraint $\theta'= 0$. 
Afterwards, in section \ref{task_dep} we remove this condition, but we fix initial and final condition on the angle $\theta$: a reaching problem can indeed require a specific direction of grasping the target object, not necessarily coincident with the one at the beginning of the movement.
For the last mentioned analysis, we will assume that the object can be reached with an arbitrary direction of the hand. Hence, from the set of geodesics connecting each couple of points, we will detect the minimimun path according to Definition \ref{definizione_distanza_p_insieme}. 
In section \ref{4.2.2}, we will consider an interval of possible directions also for the starting position. The minimum path will be modelled as the geodesic between the sets representing the conditions at the extremes (this concept is formally expressed in Definition \ref{def_geo_amm_set}).


Numerically we solve the geodesics problem 
by the use of a shooting method (see \cite{na1980computational} for further details), as follows.\\

(SHM) 
We want to solve Hamilton's equations (represented in \eqref{ham_flash} and \eqref{ham_2D}) with boundary conditions 
\begin{align}\label{bc}
\left(\gamma_i\left(0\right)\right)_{i=1,\ldots, k}=\alpha_0 
\quad,\quad\left(\gamma_i\left(T\right)\right)_{i=k+1,\ldots, n}=\alpha_1.
\end{align}
We search for a vector $\beta_0 \in\mathbb{R}^{n-k}$, which is the vector of unknown initial conditions, such that, the  $\gamma_{\beta_{0}}$ is a solution of the hamiltonian system,  \eqref{ham_flash} or \eqref{ham_2D}, with initial conditions
\begin{align}
\left(\gamma_{i, \beta_0}\left(0\right)\right)_{i=1,\ldots, n}= \left(\alpha_0, \beta_0\right)\quad \text{satisfies}\quad
\left(\gamma_{i, \beta_0}\left(T\right)\right)_{i=k+1,\ldots, n}= \alpha_1.
\end{align}
\smallskip

Finding the initial condition $\beta_0$ 
is equivalent of finding the zeros of a function of the variable $\beta_0$ 
\begin{equation}\label{fsolve}
G\left(\beta_0\right)= \left(\gamma_{i, \beta_0}\left(T\right)\right)_{i=k+1,\ldots, n}- \alpha_1.
\end{equation}

\begin{remark}
In the tests that follow, we consider an interval $\left[0, T\right]$ where we chose properly the extremum $T$ so that the functions $s\mapsto h\left(s\right)= p_{t} + v\left(s\right)p_{x}+ a\left(s\right)p_{v}\left(s\right)$ and 
$s\mapsto\psi\left(s\right)= v\left(s\right)\left(\cos\left(\theta\left(s\right)\right) p_{x}+ \sin\left(\theta\left(s\right)\right) p_{y}\right)+ a\left(s\right)p_{v}\left(s\right)+ p_{t}$ together with their derivatives were not null. Therefore, we are truly reduced to analyze geodesics of the space in accordance to their definition of admissibility, as expressed in \ref{admissible} and \ref{admissible_2D}. Moreover, in this way the solution of the initial value problem is regular, without singularities.
A different problem, which is the largest interval on which the solution are regular, has been investigated by U. Boscain et al. in \cite{boscain2012optimal} in the context of the SE(2) group.
\end{remark}

\subsection{Comparison with Flash and Hogan model}\label{comp_FH}
\begin{figure}[htbp]
\centering
\includegraphics[scale=0.35]{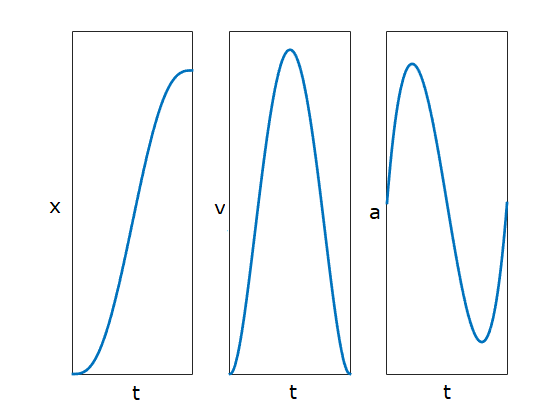}
\quad
\includegraphics[scale=0.35]{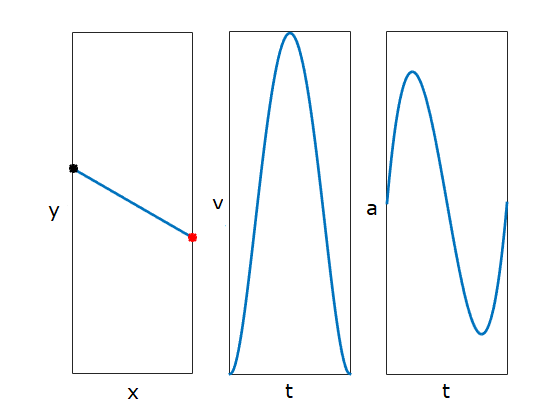}
\caption{Left. Solution's projection of system \eqref{ham_flash} over the $\left(t, x\right)$, $\left(t,v\right)$, $\left(t, a\right)$ planes, with extremes conditions \eqref{bc_2}. Right. Solution's projection of system \eqref{ham_2D} over the $\left(x, y\right)$, $\left(t,v\right)$, $\left(t, a\right)$ planes, with extremes conditions \eqref{bc3} and movement direction $\theta= \frac{5}{6}\pi$. Red and black dots respectively denote initial and final hand's position.}
\label{Flash_path}
\end{figure}

We analyze a task for which a final target is assumed to be achieved in a smooth way starting at zero speed and acceleration. Conditions at the extremes relative to velocity and acceleration match the ones analyzed by Flash and Hogan model. Hence, in the same notations of problem (SHM) and referring to geodesics equations \eqref{ham_flash}, we assume
\begin{equation}\label{bc_2}
\alpha_0= \left(t_0,x_0, v_0, a_0\right)= \left(0,x_0, 0, 0\right)\quad\text{and}\quad \alpha_1= \left(t_1, x_1, v_1, a_1\right)= \left(1, x_1, 0, 0\right).
\end{equation}
In the above conditions \eqref{bc_2}, an admissible curve exists (see Proposition \ref{p4}) and has a bell-shaped speed profile (see Remark \ref{studio_grafico}). 
In the left part of Figure \ref{Flash_path}, a representation concerning the trajectory, speed and acceleration profiles is shown.

Analogously for the two-dimensional case, by considering 
system \eqref{ham_2D} with constraint $\theta'= 0$ and boundary values
\begin{equation}\label{bc3}
\tilde{\alpha}_0= \left(t_0, x_0, y_0, \theta_0, v_0, a_0\right)= \left(0, 0, 0, \theta_0, 0, 0\right)\:,\:
\tilde{\alpha}_1= \left(t_1, x_1, y_1, v_1, a_1\right)= \left(1, x_1, y_1, 0, 0\right). 
\end{equation}
we obviously find the same curves of the predicted paths and trajectories of Flash and Hogan model 
(see Figures \ref{Flash_path} and \ref{fm} for a direct comparison). 

\begin{figure}[htbp]
\centering
\includegraphics[scale=0.4]{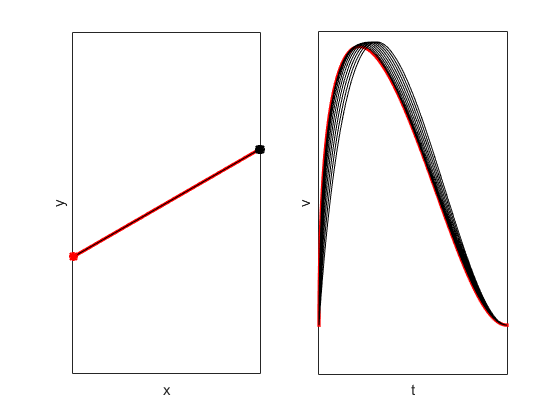}
\quad
\includegraphics[scale=0.4]{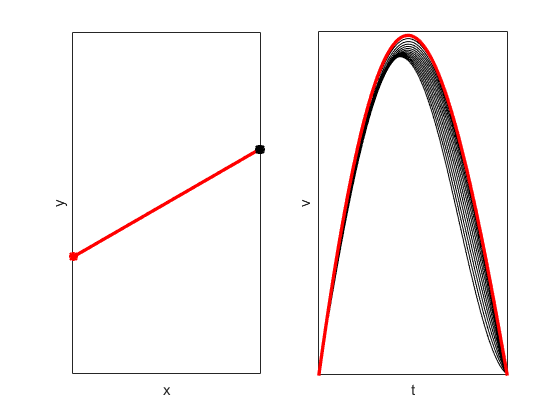}
\caption{Solutions projections referred to system \eqref{ham_2D} with constraint $\theta'= 0$ over the $\left(x,y\right)$ and $\left(t,v\right)$ planes. The red-colored path represents the geodesic with minimum length according to \eqref{def_distanza_p_insieme}. 
Left. Initial conditions are given by $\bar{\alpha}_0= \left(0, 0, 0, \frac{\pi}{6}, 0, \bar{a}_0\right)$, where $\bar{a}_0\in\left[\frac{3\pi}{10}, \frac{2\pi}{5}\right]$; final condition is assumed to be fixed at point $\tilde{\alpha}_1$, as in \eqref{bc3}. Right. Initial condition is represented by $\tilde{\alpha}_0= \left(0,0,0,\frac{\pi}{6},0, \frac{3\pi}{8}\right)$; final conditions are $\bar{\alpha}_1= \left(1,x_1, y_1, 0, \bar{a}_1\right)$, where $\bar{a}_1$ varies over the interval $\left[-\frac{2\pi}{5}, -\frac{3\pi}{10}\right]$.}
\label{a0_a1}
\end{figure}

 
Our setting allows to take into account even more general situations. For instance, we can consider a movement which does not require to start or end with a fixed acceleration. 
In the same conditions as before, we can replace $\tilde{\alpha}_0$ with  $\bar{\alpha}_0:= \left(t_0, x_0, y_0, \theta_0, v_0, \bar{a}_0\right)= \left(0, 0, 0, \theta_0, 0, \bar{a}_0\right)$, where $\bar{a}_0$ is supposed to be varying over an interval $\left[a_0^1,a_0^2\right]$, so that it is possible to analyze a range of possible initial accelerations. Then, for any choices of $\bar{a}_0$, we solve through (SHM) problem \eqref{ham_2D} with initial value $\bar{\alpha}_0$ and final value $\tilde{\alpha}_1$ and, from the set of solutions, we select the geodesic with minimum length according to distance \eqref{def_distanza_p_insieme}. In this context, the fiber of possible choices is assumed to be the one represented by the vector field $X_3$. 
In the same way, we could also think of replacing $\tilde{\alpha}_1$ with $\bar{\alpha}_1:= \left(t_1,x_1,y_1, v_1, \bar{a}_1\right)= \left(1,x_1,y_1,0,\bar{a}_1\right)$, where $\bar{a}_1\in\left[a_1^1, a_1^2\right]$. In Figure \ref{a0_a1} a representation example of these situations is shown. The red-colored curve represents the geodesic from a point to a set as in Definition \ref{definizione_distanza_p_insieme}. In this case, it clearly appears how the geodesic results to be the curve which minimizes the rate of change of acceleration over the whole speed profile. This is due to the fact that the curvature over the path $\left(x,y\right)$ is constantly null, therefore formula \eqref{length_a_2d} measures the length of a path only by taking into account the 
the tangent over the velocity function $t\mapsto v\left(t\right)$. 

\subsection{Task-dependent boundary conditions}\label{task_dep}
In this section, we generalize the class of center-out movements by exploring a set of reaching tasks which include a temporal change on the movement direction variable. We refer to subsections \ref{sub_fixed_directions} and \ref{sub_free_direction} 
for the formal arrangement of the problem.
\subsubsection{Reaching targets with prescribed directions}\label{prescribed_directions}

\begin{figure}[htbp]
\centering
\includegraphics[scale= 0.3]{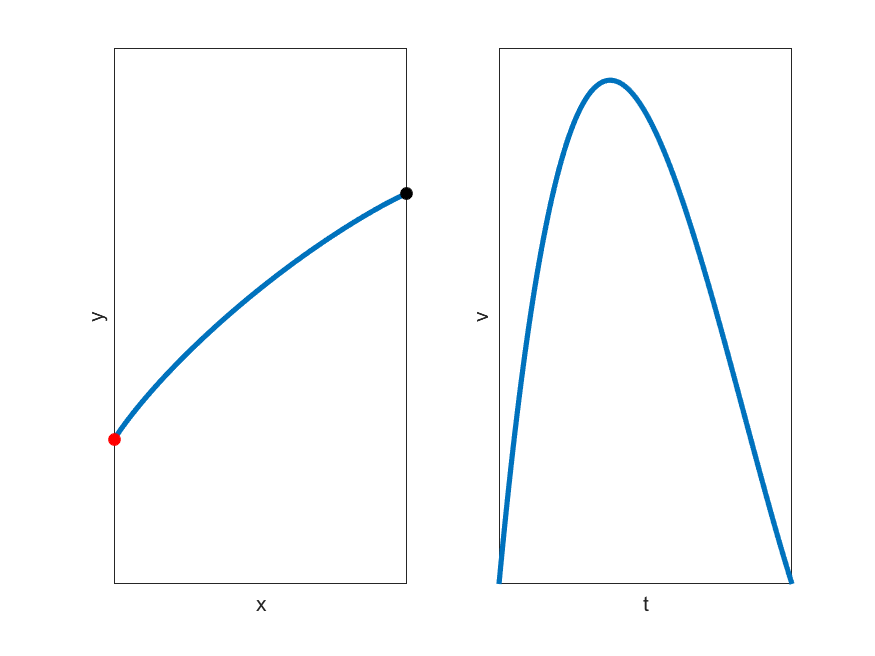}
\includegraphics[scale= 0.3]{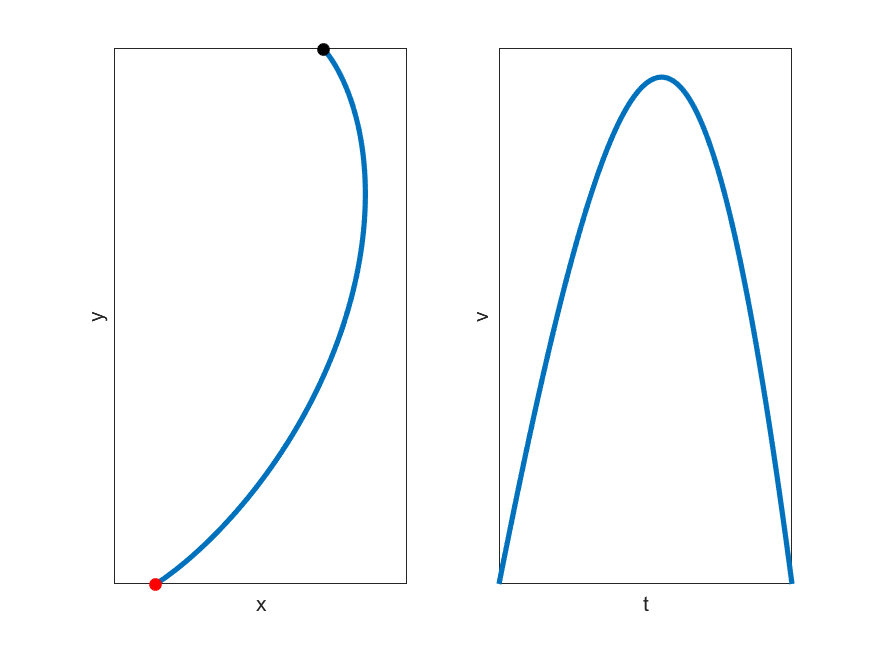}
\includegraphics[scale= 0.3]{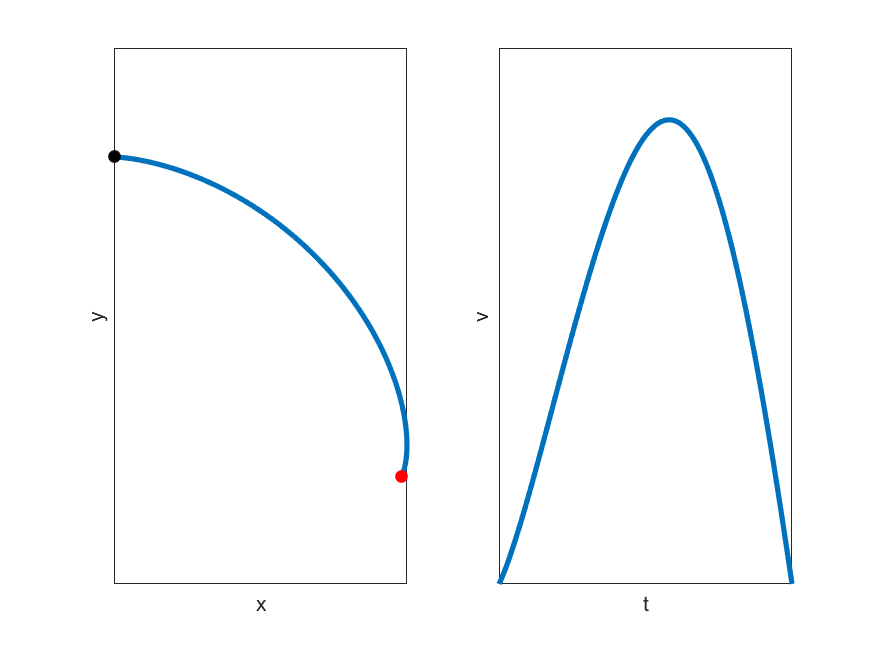}
\caption{Reaching paths and speed profiles with boundary conditions \ref{prescribed_direction}. From left to right, the assumed extreme conditions for hand's orientations $\left(\theta_0, \theta_1\right)$ are  $\left(\frac{\pi}{3}, \frac{\pi}{8}\right)$, $\left(\frac{\pi}{6}, \frac{3\pi}{4}\right)$, $\left(\frac{\pi}{4}, \pi\right)$. Accelerations couple $\left(a_0, a_1\right)$, from left to right, are $\left(\frac{\pi}{3}, -\frac{9\pi}{20}\right)$, $\left(\frac{\pi}{3}, -\frac{\pi}{4}\right)$, $\left(\frac{9\pi}{20}, -\frac{\pi}{3}\right)$.}
\label{caso3}
\end{figure}

The simulations we present aim at computing reaching movements in which the final target is assumed to be grasped with a certain orientation $\theta_1$. Moreover, an initial hand orientation $\theta_0$ is assumed to be given.  
As we did in section \ref{comp_FH}, we set up (SHM) by considering system \eqref{ham_2D} and conditions at the extremes represented by 
\begin{equation}\label{prescribed_direction}
\hat{\alpha}_0= \left(t_0, x_0, y_0, \theta_0, v_0, a_0\right)= \left(0, x_0, y_0, \theta_0, 0, a_0\right)\:,\:
\hat{\alpha}_1= \left(1, x_1, y_1, \theta_1, 0, a_1\right).
\end{equation}
Proposition \ref{prop_connection_curve_amm_2d} ensures that an admissible curve connecting $\hat{\alpha}_0, \hat{\alpha}_1$ exists, therefore we look for the missing initial conditions $\hat{\beta}_0= \left(p_{t_0}, p_{x_0},p_{y_0}, p_{\theta_0}, p_{v_0},p_{a_0}\right)$ which solve the initial value problem \eqref{ham_2D} with $\left(\hat{\alpha}_0, \hat{\beta}_0\right)$ as initial datum and which satisfy equation \eqref{fsolve}. Some examples representing the paths and velocity profiles are shown in Figure \ref{caso3}.

Even in this context, we can discuss upon the choice of initial and final accelerations. For instance, if we want to grasp an object by slowly decelerating, we loose something in terms of optimal length. Indeed, if there is a choice in terms of final accelerations in order to reach a target, the minimum path will not be the one with the smoothest features in the $\left(t,v\right)$ plane. As it is shown in Figure \ref{caso4}, if we compare those velocity profiles with increasingly steeping final accelerations, we see that those with $a_1$ more close to the $t-$axis provide smoother and longer curves. Referring to the $\left(x,y\right)$ plane, the minimum path highlighted in red is associated with a euclidean length greater than the others. This is due to the high speeds reached which, given a fixed window of times, determine a longer path to be taken. The same reasoning also apply for an interval of choices referred to initial accelerations. 


\begin{figure}[htbp]
\centering
\includegraphics[scale= 0.3]{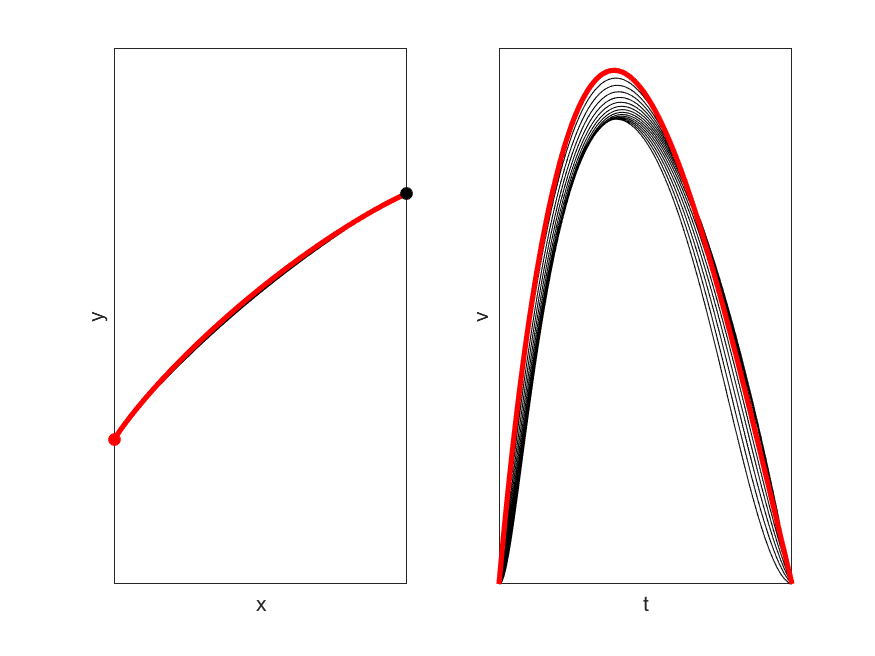}
\includegraphics[scale= 0.3]{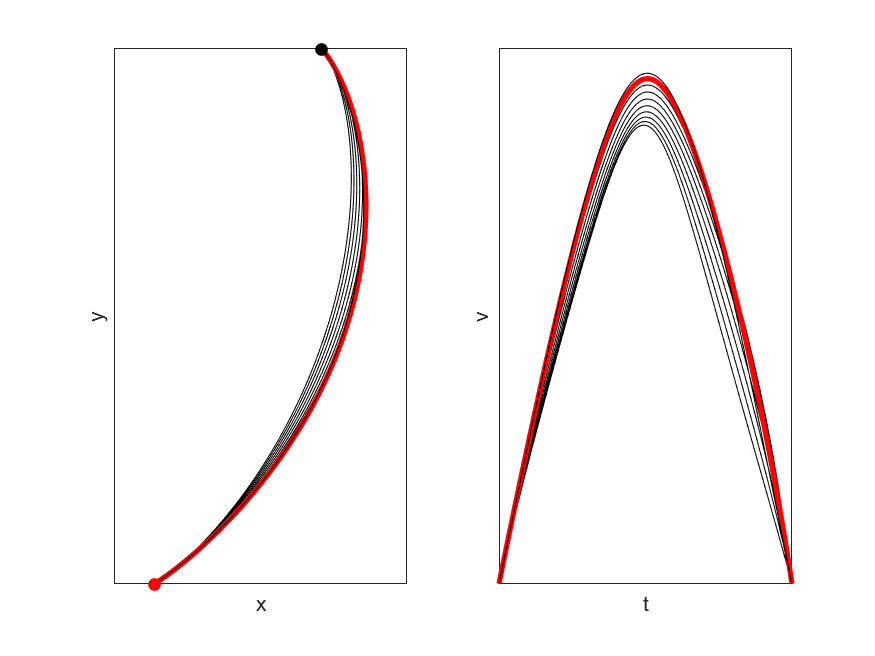}
\includegraphics[scale=0.3]{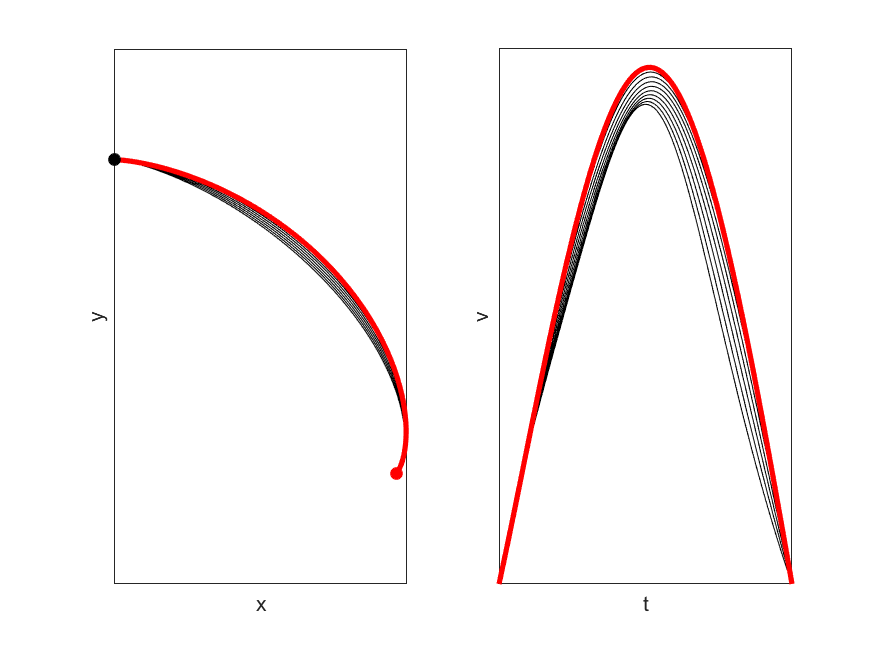}
\caption{Reaching paths and speed profiles with boundary conditions as in Figure \ref{caso3}, with final acceleration varying in intervals  $\left[-\frac{7\pi}{16}, -\frac{\pi}{3}\right]$, $\left[-\frac{\pi}{3}, -\frac{\pi}{6}\right]$, $\left[-\frac{3\pi}{8}, -\frac{\pi}{4}\right]$. The red-colored path represents the geodesic with minimum length according to \eqref{def_distanza_p_insieme}.} 
\end{figure}
\label{caso4}

\subsubsection{Reaching targets with arbitrary directions}\label{4.2.2}
In this section, we will continue the analysis by focusing on the variation of parameters with respect to the movement direction.
We point out that if we make varying both initial and final conditions it means that we are looking for a geodesic between sets, as well as if we fix an extremum and we study an interval of choice for the other boundary condition, we are searching for a geodesic from a set. We are formally considering Definitions \ref{definizione_distanza_insieme_insieme} and \ref{definizione_distanza_p_insieme} analyzed in section  \ref{sub_free_direction}. 
In a first case, we consider the cognitive situation in which there exists a range of possible movement directions in order to achieve the final target. Hence, we assume the following conditions at the extremes
\begin{equation}\label{bc4}
\hat{\alpha}_0= \left(0,x_0, y_0,\theta_0, 0,a_0\right)\quad,\quad
\hat{\alpha}^{\theta}_1= \left(1, x_1, y_1, \theta_1, 0, a_1\right),\quad\theta_1\in\left[\theta^a, \theta^b\right].
\end{equation}

\begin{figure}[htbp]
\centering
\includegraphics[scale= 0.32]{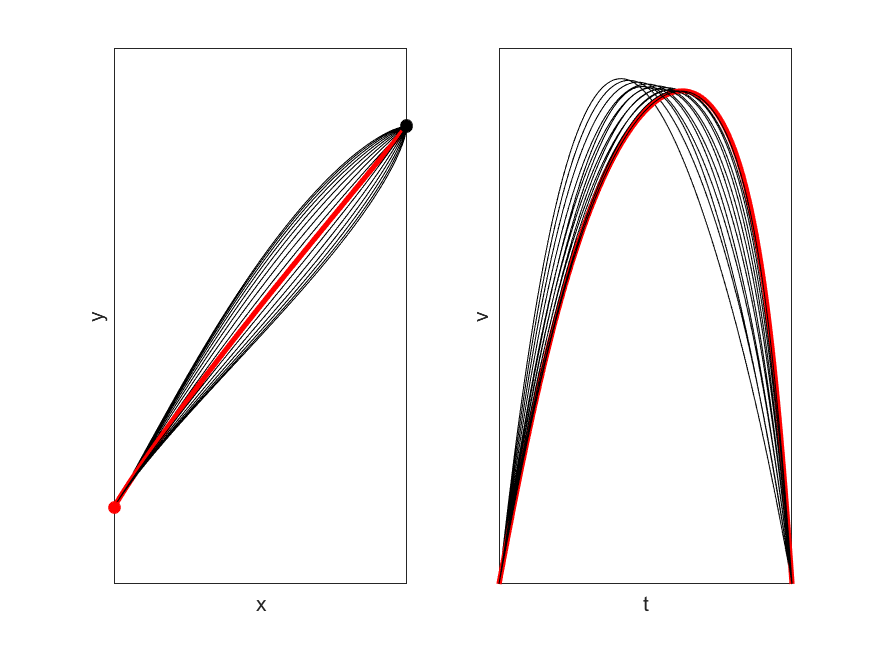}
\includegraphics[scale= 0.32]{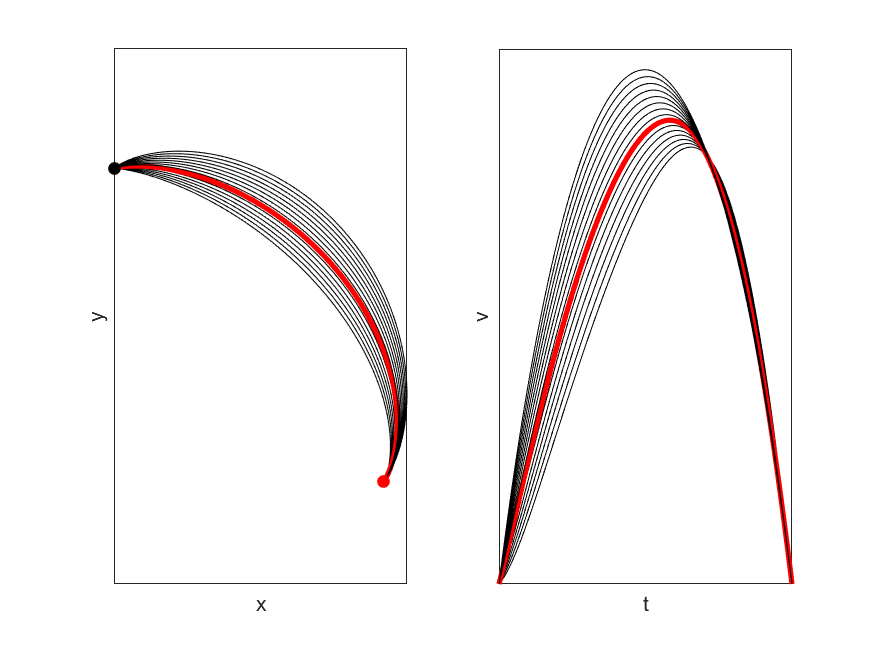}
\includegraphics[scale= 0.32]{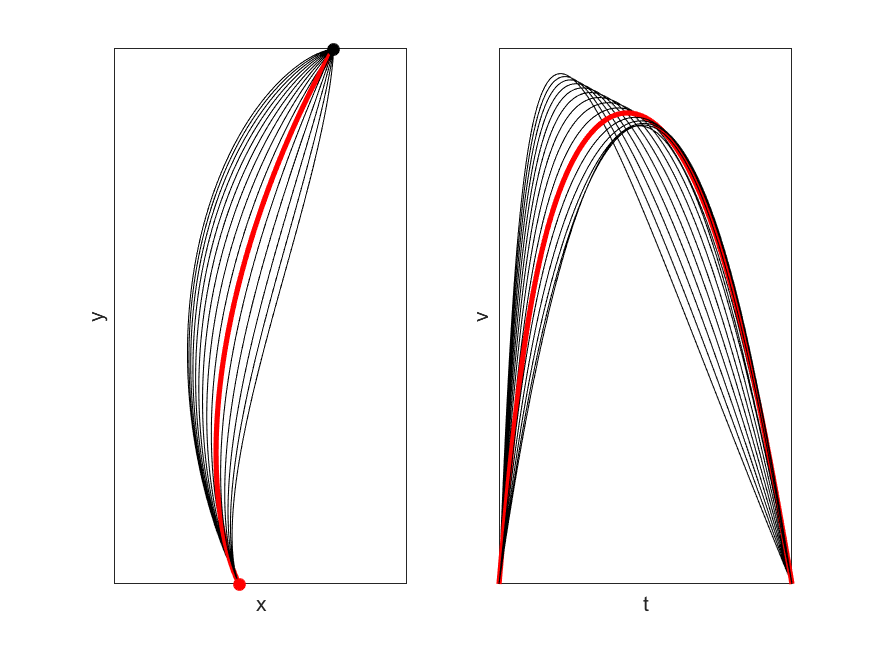}
\caption{Reaching paths and speed profiles with boundary conditions \ref{bc4}. The initial movement directions $\theta_0$ are, from left to right, $ \frac{\pi}{3},\frac{\pi}{4},\frac{2\pi}{3}$, whereas the final ones  $\theta_1$ vary on intervals $\left[0,\frac{\pi}{2}\right], \left[\pi, \frac{5\pi}{4}\right], \left[0,\frac{\pi}{2}\right]$. The assumed extreme conditions for hand's accelerations $\left(a_0, a_1\right)$ are  $\left(\frac{3\pi}{8}, -\frac{3\pi}{8}\right)$, $\left(\frac{\pi}{3}, -\frac{3\pi}{8}\right)$, $\left(\frac{\pi}{3}, -\frac{3\pi}{8}\right)$. The red colored path is the geodesic with minimum length according to \eqref{def_distanza_p_insieme}. 
}
\label{casobello}
\end{figure}

\begin{figure}[htbp]
\centering
\includegraphics[scale=0.35]{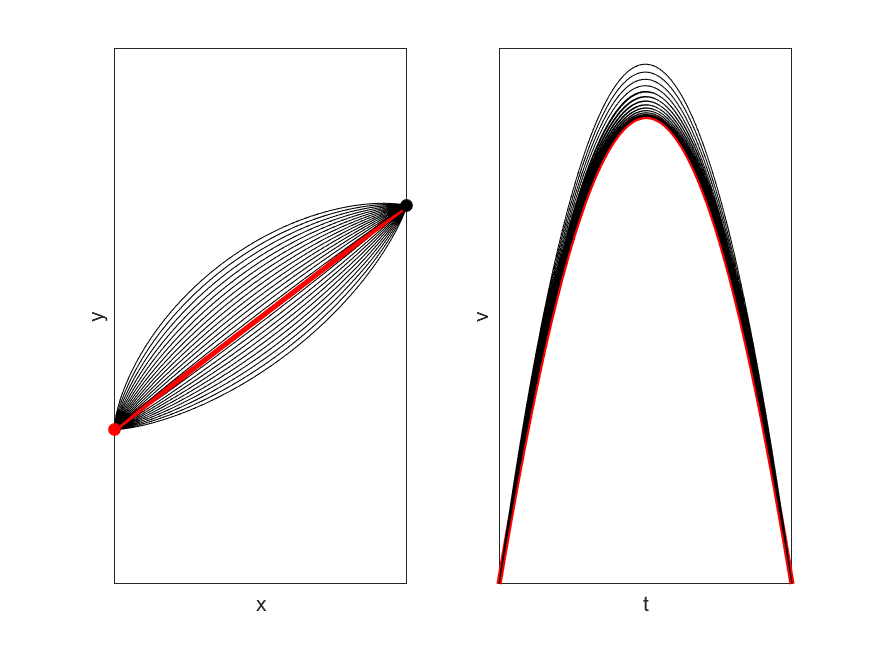}
\quad\quad
\includegraphics[scale=0.35]{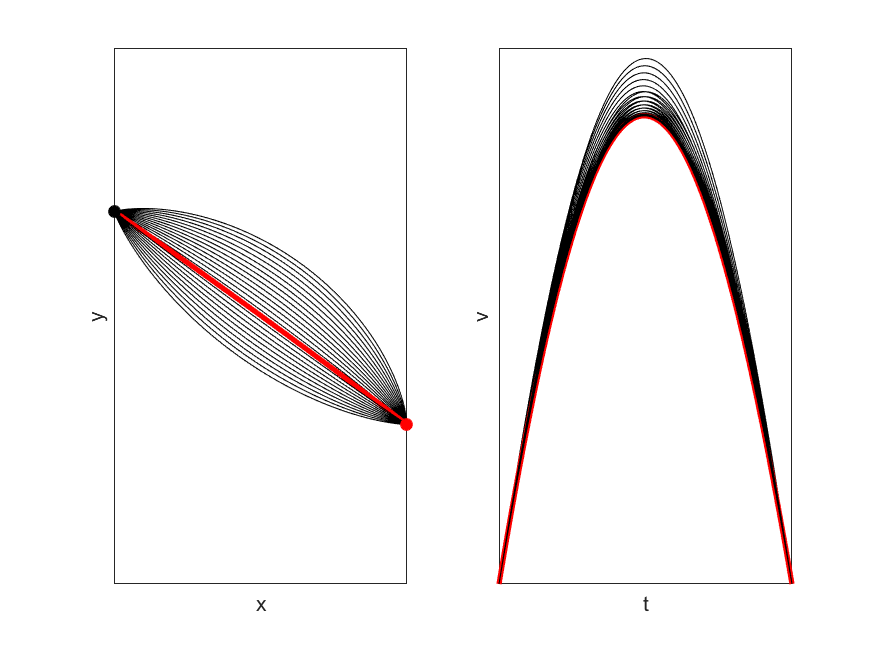}
\caption{Reaching paths and speed profiles with boundary conditions \eqref{bc5}. Initial and final hands movement directions $\theta_0, \theta_1$ are assumed to be varying in $\left[0,\frac{\pi}{2}\right], \left[-\frac{\pi}{12}, \frac{5\pi}{12}\right]$ and $\left[\frac{\pi}{2}, \pi\right], \left[\frac{11\pi}{18}, \frac{10\pi}{9}\right]$. Accelerations $\left(a_0, a_1\right)= \left(\frac{3\pi}{8}, -\frac{3\pi}{8}\right)$. 
The highlighted path is the minimum in length according to \eqref{def_distanza_insieme_insieme}.} 
\label{caso5}
\end{figure}

For any choice of $\theta_1\in\left[\theta^a, \theta^b\right]$, we solve through (SHM) the Hamiltoniam system \eqref{ham_2D} and we find out the geodesic connecting each couple of points given by $\hat{\alpha}_0$ and $\hat{\alpha}_1^\theta$. Finally, we apply \eqref{def_distanza_p_insieme} in order to catch the minimum path from the interval $\left[\theta^a, \theta^b\right]$. Some examples are represented in Figure \ref{casobello}, where the geodesics found according to Definition \ref{definizione_distanza_p_insieme} are red-colored. The geodesics account for a combination of minimum cost in terms of the spreading of curvature over the $\left(x,y\right)$ and $\left(t,v\right)$ planes, as it is inherited from Definition \ref{length_a_2d}. As we can expect, if we assume a freedom of choice even for the initial movement directions, as in the following boundary conditions
\begin{equation}\label{bc5}
\hat{\alpha}_0^\theta= \left(0, x_0, y_0,\theta_0, 0, a_0\right)\:,\:
\hat{\alpha}_1^\theta= \left(1,x_1, y_1, \theta_1, 0, a_1\right),\quad\theta_0\in\left[\theta^a, \theta^b\right],\,\theta_1\in\left[\theta^c,\theta^d\right],
\end{equation}
the resulting geodesics turn out to be curves near to be straight paths which account for the minimum difference between the  starting and ending directions, as it is shown in Figure \ref{caso5}.

\section{Conclusions}\label{concl}
Our model takes Flash and Hogan's phenomenological model as a point of reference. Their model is based on the selection of a cost function whose minima are found to be in good agreement with the simplest and smoothest motion trajectories. 
We recovered the minimizing trajectories found by Flash and Hogan as geodesics in an Engel-type manifold. Specifically we consider only a specific class of curves, called admissible, which are lifting of mono-dimensional paths. For these curves, we provided a connectivity property and the existence of length minimizers.
Admissible geodesics represent our model for center-out type movements. 
We then extended the previous model 
by considering admissible curves into a new non-nilpotent subriemannian structure. In this second part, we proved the same results of connectivity and we further studied their regularity. In section \ref{results} we showed a qualitative analysis on how admissible geodesics allow to recover a broader variety of reaching task. 

\newpage
\addcontentsline{toc}{section}{References}
\bibliographystyle{plain}
\bibliography{neurogeometry_reaching_bib}

\end{document}